\documentclass[12pt, a4paper]{amsart}
\usepackage{amsmath,amssymb,amscd,amsfonts}

\newtheorem{thm}{Theorem}[section]
\newtheorem{lem}[thm]{Lemma}
\newtheorem{rem}[thm]{Remark}
\newtheorem{prop}[thm]{Proposition}
\newtheorem{cor}[thm]{Corollary}

\newcommand{\Pic}{{\text{\rm Pic}}}
\newcommand{\PX}{{\text {\rm Pic}}^0(X)}
\newcommand{\DI}{{\text {\rm dim}}   }

\newcommand{\DX}{{\text {\rm dim}}(X)}
\newcommand{\DY}{{\text {\rm dim}}(Y)}

\newcommand{\alb}{\text{\rm alb}}
\newcommand{\Alb}{\text{\rm Alb}}
\newcommand{\lra}{\longrightarrow}
\newcommand{\ot}{{\otimes}}
\newcommand{\OO}{{\mathcal O}}
\newcommand{\OY}{{\mathcal{O} _Y}}
\newcommand{\OZ}{{\mathcal{O} _Z}}
\newcommand{\OX}{{\mathcal{O} _X}}
\newcommand{\dual}{^{\vee}}
\newcommand{\inv}{^{-1}}

\newcommand{\epsi}{\epsilon}
\newcommand{\si}{\sigma}
\newcommand{\Si}{\Sigma}

\newcommand{\ddual}{^{\vee\vee}}
\newcommand{\F}{{\mathcal F}}
\newcommand{\PP}{{\mathcal{P}}}
\newcommand{\PPP}{{\mathbb{P}}}

\newcommand{\ox}{{\omega _X}}
\newcommand{\FF}{{\mathcal{F}}}
\newcommand{\GG}{{\mathcal{G}}}
\newcommand{\KK}{{\mathcal{K}}}
\newcommand{\QQQ}{{\mathcal{Q}}}
\newcommand{\QQ}{{\mathbb{Q}}}

\newcommand{\V}{\mathcal V}
\newcommand{\sh}{{\hat {\mathcal S}}}

\title[On the birational geometry ...]{On the birational geometry 
of varieties of maximal
Albanese dimension}
\author{Christopher D. Hacon \and Rita Pardini}
\date{}

\begin{document}
\maketitle

\section{Introduction}

Combining the generic vanishing theorems of Green and Lazarsfeld \cite{GL1},
\cite{GL2},
the theory of Fourier Mukai transforms \cite{Mu} and the results of Koll\`ar
on higher direct images of dualizing sheaves \cite{Ko1}, \cite {Ko2},
it is possible to obtain surprisingly precise information about
the birational geometry of varieties of maximal Albanese dimension.
In \cite{CH1} and \cite{CH3} for example, it is shown that the following
conjecture of Koll\`ar holds (cf. \cite{Ko3})

\medskip
\noindent{\bf Theorem.} {\em If $X$ is a smooth complex algebraic variety with
$P_2(X)=1$ and $q(X)=\dim (X)$, then $X$ is birational to an abelian variety.}
\medskip

In this paper we show how these techniques can be used to give
effective criteria for a morphism of varieties of maximal Albanese dimension
to be birational.
In particular we prove the following

\medskip
\noindent{\bf Theorem 1.} {\em Let $f\colon X\to Y$ be a generically finite
morphism of smooth complex projective
$n$-dimensional varieties of maximal Albanese dimension.
If the induced maps
$$H^i(Y,\omega _Y\ot P)\to H^i(Y,f_*\omega _X\ot P)$$
are isomorphisms for all $i\geq 0$ and all $P\in \Pic ^0(Y)$,
then $X\to Y$ is a birational morphism.}
\medskip

As an application of this, we prove a refinement of \cite{Hac} Theorem 3.1.
\medskip

\noindent{\bf Corollary 2.} {\em Let $\nu\colon X\to A$
be a morphism from a smooth variety
of dimension $n$ to an abelian variety of dimension $n+1$ such that
$h^0(X,\omega _X \ot \nu ^*P)=1$ for all $\OO _A \ne P\in \Pic ^0(A)$
and the map $\nu ^*\colon H^0(A,\Omega ^{n}_A)\to H^0(X,\Omega ^{n}_X)$ is an
isomorphism. Then $\nu$ is birational onto its image $\bar {X}$ and
$\bar {X}$ is a principal polarization.}
\medskip

It is well known that if $f\colon X\to Y$ is a morphism of $n$-dimensional
smooth projective varieties of general type such that $P_m(X)=P_m(Y)$ for
$m\gg 0$ then
$f$ is birational. If $Y$ is of maximal Albanese dimension, it suffices
to verify the above condition for $m=2$.
\medskip

\noindent{\bf Theorem 3.} {\em Let $f\colon X\to Y$ be a dominant
morphism of smooth $n$-dimensional complex projective
algebraic varieties, $Y$ of general type and of maximal Albanese dimension.
If $P_2(X)=P_2(Y)$, then
$f$ is a birational morphism.}
\medskip

Next we turn our attention to the Albanese morphism. By a
theorem of Kawamata \cite{Ka}, if $\kappa (X)=0$ then the Albanese map of $X$
is surjective. In particular $\dim(X)\geq q(X)$. If furthermore
$q(X)=\dim (X)$ then $\alb _X\colon X\to \Alb (X)$ is a birational morphism.
In \cite {Ko4}, Koll\`ar proves effective versions of this result.
We further refine \cite[Theorem 11.2]{Ko4}
\medskip

\noindent{\bf Theorem 4.}
{\em Let $X$ be a smooth projective variety. If $P_2(X)=1$ or
$0<P_m(X) \leq 2m-3$ for some $m\geq 3$, then
$\alb _X\colon X\to \Alb (X)$ is surjective. }
\medskip

As remarked above, the case $P_2(X)=1$, $q(X)=\dim (X)$
was completely understood
in \cite{CH1} and \cite{CH3}. We study the case $P_3(X)=2$, $q(X)=\dim (X)$.
\medskip

\noindent{\bf Theorem 5.} {\em Let $X$ be a smooth projective variety.
\noindent Then
$P_3(X)=2$ and $q(X)=\dim (X)$ iff:
\begin{itemize}
\item[a)] there is a surjective map $\phi\colon \Alb(X)\to E$, where $E$ 
is a curve of
genus
$1$;

\item[b)]$\alb _X\colon X\to\Alb (X)$
is birational to a smooth double cover of $\Alb (X)$ defined by
$Spec (\OO _{\Alb (X)}\oplus \phi^*L\inv\ot P)$,
where $L$ is a line bundle of $E$ of
degree $1$ and
$P\in
\Pic ^0(X)\setminus \phi^*\Pic^0(E)$ is
$2$-torsion. The branch locus of the cover is the union of two
distinct fibers of $\phi$.
\end{itemize}}
\medskip

\noindent{\bf Acknowledgments.} This research was started in June 2000,
during a visit of
the first author to Pisa, supported by G.N.S.A.G.A. of C.N.R.
The first author is grateful to A. J. Chen and R. Lazarsfeld for
helpful conversations. We are also in debt 
to Ein and Lazarsfeld for \cite{EL2} and their permission to include Lemma D,
\cite{EL2} in this paper.
\medskip

\noindent{\bf Notation and conventions.}
We work over the field of complex numbers.
We identify Cartier divisors and line bundles on
a smooth variety, and we use the additive and
multiplicative notation interchangeably. If
$X$ is a smooth projective variety, we let $K_X$ be a  canonical divisor,
so that
$\ox=\OO_X(K_X)$, and we denote by
$\kappa(X)$ the Kodaira dimension,  by
$q(X):=h^1(\OO_X)$ the {\em irregularity} and by $P_m(X):=h^0(\omega _X^m)$
the {\em $m-$th
plurigenus}. We denote by  $\alb_X\colon X\to \Alb(X)$ the Albanese map  and by
$\Pic^{\tau}(X)$ the subgroup of torsion elements of $\Pic^0(X)$.
For a $\QQ-$divisor $D$
we let $\lfloor D\rfloor$ be the integral part and
$\{D\}$ the fractional part. Numerical
equivalence is denoted by $\equiv$ and we write $D\prec E$ if
$E-D$ is an effective divisor. If $f\colon X\to Y$ is a
morphism, we write $K_{X/Y}:=K_X-f^*K_Y$ and we  often denote by  $F_{X/Y}$
the general fiber of
$f$.  The rest of the notation is
standard in algebraic geometry.

\section{Preliminaries}
Here we recall several results from the literature and we prove some technical
statements that will be needed later.

\subsection{The Albanese map and the Iitaka fibration}
Let $X$ be a smooth projective variety. If $\kappa(X)>0$,
then the Iitaka fibration of $X$
is a morphism of projective varieties $f\colon X'\to Y$, with $X'$ 
birational to $X$ and
$Y$  of dimension
$\kappa(X)$, such that the general fiber of $f$ is smooth,
irreducible, of zero
Kodaira dimension.  The Iitaka fibration is determined only up to
birational equivalence.
Since we are interested in questions of a birational nature,
we usually assume that
$X=X'$ and that $Y$ is smooth.

$X$ has {\em  maximal Albanese dimension} if
$\DI {(\alb _X(X))}=\DI (X)$.
\begin{prop}\label{albanese}Let $X$ be a smooth projective variety of maximal
Albanese dimension, and let
$f\colon X\to Y$ be the Iitaka fibration (assume $Y$ smooth). Denote by 
$f_*\colon \Alb
(X)\to \Alb (Y)$ the homomorphism induced by
$f$ and  consider the commutative diagram:
$$
\CD
X &@>{\alb _X}>> & \Alb (X)\\
@V{f}VV& & @V{f_*}VV \\
Y &@>{\alb _Y}>> & \Alb (Y).\\
\endCD
$$

Then:
\begin{itemize}

\item[a)] $Y$ has maximal Albanese dimension;

\item[b)] $f_*$ is surjective and $\ker f_*$ is connected of  dimension
$\dim (X)-\kappa(X)$;

\item[c)] there exists an abelian variety
$P$ isogenous to $\ker f_*$ such that the general fiber of $f$ is birational to
$P$.
\end{itemize}
\end{prop}
\begin{proof}
The dual of the differential of $f_*$ at $0$ is the pull--back map
$H^0(Y,\Omega^1_Y)\to H^0(X,\Omega^1_X)$, which
is clearly injective. It follows
that $f_*$ is surjective. Denote by $K$ the connected component of $\ker f_*$
that contains $0$ and set
$A:= \Alb (X)/K$. The abelian variety $A$  is a finite
\'etale cover of $\Alb (Y)$, and
$f\colon X\to Y$ factors through the induced map $Y\times_{
\Alb (Y)}A\to Y$, which
is a connected \'etale cover of the same degree as $A\to \Alb
(Y)$. Since the fibers
of
$f$ are connected by assumption, it follows that $A\to \Alb (Y)$ is an 
isomorphism
and the fibers of $f_*$ are connected.

Let $F$ be a general fiber of $f$. By Theorem $1$ of
\cite{Ka}, the Albanese map of $F$ is surjective, and therefore
the image of
$F$ via $\alb_X$ is a translate of an abelian subvariety of $\Alb (X)$ which
we again denote by $K$.
$K$ is
contained in  $\ker f_*$ and is independent of $F$, since $F$ moves
in a continuous system.
Denote by $A$ the quotient abelian variety $\Alb(X)/K$.
The induced map $X\to A$ is
constant on the general fiber of $f$ and thus induces a rational
map $\phi\colon Y\to
A$. Since $A$ is an abelian variety, $\phi$ is actually a morphism
and $\phi(Y)$
generates $A$ by construction. It follows that $\dim (A)\le q(Y)=\dim
\Alb(X)-\dim(\ker f_*)$, namely $\dim (K)\ge
\dim(\ker f_*)$. It follows that $K$ is equal to
$\ker f_*$.
By the theorem on the dimension of the fibers of a morphism, if $X$ has
maximal Albanese dimension, then $Y$ also has maximal Albanese dimension, and
$\ker f_*$  has dimension $\dim (X)-\dim (Y)=\dim (X)-\kappa(X)$.
Thus $q(F)\ge \dim (F)$,  and
Corollary
$2$ of
\cite{Ka} implies that  the Albanese map of
$F$ is a birational morphism. So the Albanese variety of $F$ is isogenous
to $\ker
f_*$ and, in particular, it does not depend on $F$, since $F$ moves in a
continuous system.
 \end{proof}

\subsection{Fourier Mukai Transforms}\label{mukai}

Let $A$ be an abelian variety, and denote the corresponding dual
abelian variety by $\hat {A}$.
Let $\PP$ be the normalized Poincar\'{e} bundle on
$A\times \hat {A}$ and let $\pi _{\hat {A}}\colon A\times \hat{A}\to \hat{A}$ 
be the
projection. For any point
$y\in
\hat {A}$, let
$\PP _y$ denote the associated topologically trivial line bundle.
Define the functor $\sh $ of $\OO_A$-modules
into the category of $\OO_{\hat {A}}$-modules by
$$\sh (M)=\pi _{\hat {A},*}(\PP \ot \pi ^*_AM).$$
The derived functor $R\sh$ of $\sh$ then induces an equivalence of
categories between the two derived categories $D(A) $ and $D(\hat {A})$.
In fact, if $R{\mathcal S}$ is the analogous functor on $D(\hat{A})$, one has 
(cf. \cite{Mu}): {\it There are isomorphisms of functors:
$$R{\mathcal S}\circ R\sh \cong (-1_A)^*[-g]$$
and
$$R\sh \circ R{\mathcal S} \cong (-1_{\hat {A}})^*[-g],$$
where $[-g]$ denotes ``shift the complex $g$ places to the right''.}

The {\em index theorem} (I.T.) is said to hold for 
a coherent sheaf $\FF$ on $A$
if there exists an integer $i(\FF )$ such that for all $j
\ne i(\FF )$,
$H^j(A,\FF \ot P)=0$ for all $P\in \Pic ^0(A)$.
The {\em weak index theorem} (W.I.T.) holds for a coherent sheaf $\FF$ if
there exists an integer, which we again denote by $i(\FF )$, such that 
for all $j \ne i(\FF )$, $R^j\sh
(\FF )=0$. It is easily seen that the I.T. implies the W.I.T.
If $\FF$ satisfies the W.I.T., we  denote the coherent
sheaf $R^{i(\FF )}\sh (\FF )$ on $\hat {A} $ by $\hat {\FF }$. 
By Corollary 2.4 of
\cite{Mu}, $\hat{\FF}$ also satisfies W.I.T. and has index $\dim(A)-i(\FF)$.

Here are some  consequences of the theory of \cite{Mu}:
\begin{lem}\label{nonvanish}{(Non-vanishing)}
If $\FF $ is a coherent sheaf on $A$ such that $h^i(\FF \ot P)$ $=0$ for all $i
\geq 0$ and all $P\in \Pic ^0(A)$, then $\FF =0$.
\end{lem}
\begin{proof} Follows immediately from \cite[Theorem 2.2] {Mu}.
\end{proof}
\begin{prop}\label{inclusion} Let $\psi\colon \FF \hookrightarrow
\GG$ be an inclusion of coherent sheaves on $A$
inducing isomorphisms $H^i(A,\FF \ot P)\to  H^i(A,\GG \ot P)$ for all
$i\geq0$ and all $P\in \Pic ^0(A)$. Then $\psi$ is an isomorphism of sheaves.
\end{prop}
\begin{proof} Let $\KK$ be the cokernel of $\psi\colon\FF \hookrightarrow
\GG$. Clearly $h^i(A,\KK \ot P)=0$ for all $i
\geq 0$ and all $P\in \Pic ^0(A)$. Therefore, $\KK =0$ by 
Lemma \ref{nonvanish}, and
$\psi$ is an isomorphism.
\end{proof}
\begin{prop}\label{abelsupport}
Let  $\FF $ be a coherent sheaf on $A$ such that for all $P\in \Pic ^0(A)$
we have $h^0(\FF \ot P)=1$ and $h^i(\FF \ot P)=0$ for all $i>0$.
Then $\FF $ is supported on an abelian subvariety of $A$.
\end{prop}
\begin{proof}
$\FF$
satisfies the I.T., $M:=\hat {\FF\  }$ is a
line bundle  on $\Pic ^0 (A)=\hat {A}$
that satisfies the  W.I.T. with index $i(M)=\dim (A)$,
and $\hat {M }=(-1_A)^*\FF$.
Any line bundle with $i(M)=\dim (A)$ is negative semidefinite.
It is well known that there exists a morphism of abelian varieties
$b\colon \hat{A}\to A'$ such that $M=b^*M'$ for some negative definite
line bundle $M'$ on $A'$.
It follows that $\hat {M }$ and hence $\FF$ are supported on
the image of $b^*\colon \hat{A'}\to A$.
\end{proof}
We also recall the following result from \cite{Hac}:
\begin{thm}\label{theta}
Let $X$ be a smooth complex projective variety of dimension $n$,
let $A$ be an abelian
variety of dimension $n+1$, and let
$f\colon X\to A$ be a morphism generically finite onto its image. Assume that:
\begin{itemize}
\item[a)] $h^0(X, \Omega^i_X)=$ $\binom{n+1}{i}$ for $0\le i\le n$;

\item[b)] $h^i(X,
\omega_X\ot(f^*P))=0$ for all $\OO_A\ne P\in \Pic^0(A)$ and all $i>0$.
\end{itemize}
\noindent Then $A$ is
principally polarized and $f(X)$ is a theta divisor.
\end{thm}
\begin{proof}
Cf. \cite{Hac}, Corollary 3.4.
\end{proof}
\subsection{Direct images of dualizing sheaves}\label{dualizing}

The following Theorem summarizes the most important facts  
on direct images of adjoint
bundles. We refer to  \cite{Ko3}, (10.1.5) for the definition 
of ``klt'' (Kawamata log
terminal).

\begin{thm} \label{omegavan}Let $f\colon X\to  Y$ be a  surjective 
map of projective
varieties,
$X$ smooth and
$Y$ normal.
Let $M$ be a line bundle on $X$
such that $M\equiv f^* L +\Delta$, where $L$
is a $\QQ-$divisor on $Y$
and $(X,\Delta )$ is klt. Then
\begin{itemize}
\item[a)] $R^jf_*(\omega _X\ot M) $ is torsion free for $j\geq 0$;

\item[b)] If $L$ is nef and big, then $H^i(Y,R^jf_*(\omega _X\ot M))
=0$ for $i>0$, $j\geq 0$;

\item[c)] $R^if_* \omega _X$ is zero for $i>\DX -\DY$.

\end{itemize}
\end{thm}
\begin{proof} Statements a)  and b)  correspond to Corollary  
10.15 of \cite{Ko3}.
Statement c) follows from Theorem 2.1 of \cite{Ko1}.
\end{proof}
The following observation is an easy consequence of the previous results.
\begin{lem}\label{Ri}Let $f\colon X\to Y$ be a surjective map of projective 
varieties, with
$X$ smooth and $Y$ normal,  and let
$M$ be a line bundle on $X$ such $M\equiv f^*L+\Delta$ , 
where $L$ is a nef e big
$\QQ-$divisor on $Y$ and $(X,\Delta)$ is klt.
If  $g\colon Y\to Z$ is any   morphism with $Z$ projective, then

$$H^i(Y,R^jf_*(K_X+M))=H^i(Z, g_*(R^jf_*(K_X+M)) )\quad i,j\ge 0.$$
\end{lem}
\begin{proof}
Clearly we may replace $Z$ by the image of $g$ and assume that $g$ is onto. 
In addition,  may assume that $Z$ is normal, since $g$ factors through the
normalization $\nu\colon Z'\to Z$ and $R^i\nu_*\FF=0$ for any 
coherent sheaf $\FF$ on
$Z'$ and $i>0$. By the Leray spectral sequence, to prove the claim 
it is enough  to
show that  the sheaves
$R^ig_*(R^jf_*(\omega_X\ot M))$ are zero for $j\ge 0$ and   
$i>0$. Fix an ample divisor
$H$ on
$Z$ such that: 

a)
$R^ig_*(R^jf_*(\omega_X\ot M))\ot H$ is generated by global sections; 

b)
$H^k(Z, R^ig_*(R^jf_*(\omega_X\ot M))\ot H)=0$ for $k>0$.

 The Leray spectral sequence that computes
$H^i(Y, R^jf_*(\omega_X\ot M)\ot g^*H)$ degenerates at
$E_2$ by condition b), and we have isomorphisms between 
$H^i(Y, R^jf_*(\omega_X\ot
M)\ot g^*H)$ and $H^0(Z, R^ig_*(R^jf_*(\omega_X\ot M))\ot H)$. For $i>0$, the
former group vanishes by Theorem \ref{omegavan}, b), since $L\ot g^*H$ 
is again nef and big. Now for $i>0$ the   sheaf
$R^ig_*(R^jf_*(\omega_X\ot M))\ot H$   is zero  by condition a), and thus
$R^ig_*(R^jf_*(\omega_X\ot M))$ is also zero.
\end{proof}
We use the techniques of \ref{mukai}  to reprove
a particular case of a result of Koll\`ar,
(\cite{Ko3}, Theorem 14.7). We remark that in
\cite{Ko3}
 the hypothesis that  $Y$ be  of
maximal Albanese dimension is replaced by the less restrictive hypothesis
that $Y$ has generically large fundamental group. However
the proof is more technical.
\begin{thm}\label{kollar}{\em (Koll\`ar)}
Let $f\colon X\to  Y$ be a
surjective map of projective varieties, $X$ smooth and $Y$ normal.
Let $M$ be a line bundle on $X$
such that $M\equiv f^* L +\Delta$, where $L$
is a nef and big $\QQ$-divisor on $Y$
and $(X,\Delta )$ is klt. If $Y$ is of maximal Albanese dimension,
then for every $j\ge 0$:
\begin{itemize}

\item [a)] $h^0(Y,R^jf_*(\omega
_X\ot M)\ot P)$  is independent of $P\in \Pic^0(Y)$;

\item [b)] $H^0(Y,R^jf_*(\omega
_X\ot M)\ot P)=0$
iff $R^jf_*(\omega _X\ot M      )=0$.
\end{itemize}
\end{thm}

\begin{proof}
If
$R^jf_*(\omega _X\ot M)=0$, then clearly $H^0(Y,R^jf_*(\omega
_X\ot M)\ot P)=0$ for all $P\in \Pic^0(Y)$.

Suppose now that
$R^jf_*(\omega _X\ot M)\ne 0$. The sheaf
$R^jf_*(\omega _X\ot M)$ is torsion
free by Theorem \ref{omegavan}, a). Let $a\colon   Y
\to A:=\Alb(Y)$ be the Albanese map of
$Y$, which is generically finite by assumption.   The sheaf
$a_*(R^jf_*(\omega _X\ot M ))$ is   non-zero,
supported on $a(Y)$.
By Lemma \ref{Ri}, for every $P\in \Pic^0(Y)$ and for every  $i\ge 0$ we have
isomorphisms
 $H^i(A,a_*(R^jf_*\omega_X\ot M)\ot
P)\cong H^i(Y, R^jf_*(\omega_X\ot M)\ot P)$.
If $i>0$, then  the latter group vanishes
by Theorem \ref{omegavan}, b),  and hence
$$h^0(Y, R^jf_*(\omega_X\ot M)\ot
P)=\chi(R^jf_*(\omega_X\ot M)\ot P)=
\chi(R^jf_*(\omega_X\ot M))$$ is independent of
$P$. In addition,  by Lemma
\ref{nonvanish} there is
$P\in \Pic ^0(Y)$ such that
$h^0(A,a_*(R^jf_*(\omega _X\ot M ))\ot
P)>0$. So
$h^0(Y,R^jf_*(\omega _X\ot M )\ot P)$ $=h^0(A,a_*(R^jf_*(\omega_X\ot
M)\ot P))>0$   for every
$P\in \Pic^0(Y)$.
\end{proof}

We  also need to understand the behavior of direct
images of dualizing sheaves with respect to the Stein factorization.

Let $Y$ be a smooth variety. For a coherent sheaf $\F$ on $Y$, we denote by
$\F ^{\dual}$ the {\it dual sheaf\ } $Hom_{\OO_Y}(\F,\OY)$. A sheaf $\F$ is
said to be
{\it reflexive} if the natural map $\F\to\F\ddual$ is an isomorphism.
Following \cite{OSS}, we say that $\F$ is {\em normal} iff for every open set
$U\subset Y$ and every closed subset $C\subset Y$ of codimension $>1$ the
restriction map
$\F(U)\to \F(U\setminus C)$ is an isomorphism.
\begin{prop}\label{OXfree} Let $X$, $Y$ be smooth projective varieties,
$f \colon X\to Y$  a
generically finite morphism, and  let
$X\stackrel{h}\to Z\stackrel{g}\to Y$ be the
Stein factorization of $f$.
If $f_*\omega_X$ is locally free on $Y$, then $g$ is a
flat morphism and $g_*\OZ=(f_*\omega_X)\dual\otimes\omega_Y$.
\end{prop}
\begin{proof} By the definition of Stein factorization,
$h$ is birational,  $g$ is  finite
and $g_*\OZ=f_*\OX$. The variety $Z$ is normal, since $X$ is normal, hence
$g_*\OZ$ is a normal sheaf, as defined above. Denote by $\omega_Z$ the
dualizing sheaf of
$Z$. By
\cite{Re}, Corollary (8) p.
283, if
$j\colon Z_0\to Z$ is the inclusion of the smooth locus $Z_0$ of $Z$, one has
$\omega_Z=j_*\omega_{Z_0}$. If we write $X_0=h\inv(Z_0)$ and denote by
$h_0\colon
X_0\to Z_0$ the restriction of $h$ to $Z_0$, then we have
$(h_0)_*\omega_{X_0}=\omega_{Z_0}$ since $h_0$ is a birational morphism
of smooth
varieties. This identification induces a sheaf map $\psi\colon f_*\omega_X\to
g_*\omega_Z$ that is an isomorphism outside the image $\Si$ of the
singular locus of
$Z$, which is a closed subset of $Y$ of codimension $>1$. 
The sheaf $f_*\omega_X$ is
normal, since it is locally free and $Y$ is smooth, 
while $g_*\omega_Z$ is normal by
the definition of $\omega_Z$ and by the fact that $g$ is finite.    
It follows that  $\psi$
is an isomorphism.

By \cite{Hart}, Chapter III,
ex. 6.10b, there is an isomorphism
$g_*(g^!\omega_Y)\cong
(g_*\OZ) \dual\otimes\omega_Y$, where, for a coherent sheaf $\F$ of $Y$ and a
finite morphism $g\colon Z\to Y$, one defines
$g^!\F$ as the coherent sheaf on $Z$ corresponding to the
$g_*\OZ-$module $Hom_{\OY}(g_*\OZ,\F)$. Now
$g^!\omega_Y=\omega_Z$, by \cite{Hart}, Chapter
III, ex. 7.2, hence $g_*\omega_Z=(g_*\OZ) \dual\otimes\omega_Y$. 
Taking duals, one
has $(g_*\OZ)\ddual\cong
(g_*\omega_Z)\dual\otimes\omega_Y\cong(f_*\omega_X)\dual\otimes\omega_Y
$.
The sheaf  $g_*\OZ$ is normal and torsion free, since $Z$ is normal, 
and thus it is
reflexive  by Lemma 1.1.12 of \cite{OSS}. Thus we have
$g_*\OZ\cong (f_*\omega_X)\dual\otimes\omega_Y$. In particular, $g_*\OZ$ is
locally free, namely $g$ is flat.\end{proof}

\subsection{Cohomological Support Loci}\label{supportloci}

Let $\pi\colon X\to A$ be a morphism from a smooth projective variety $X$ to
an abelian variety $A$.
If $\FF$ is a coherent sheaf on $X$,
then one can define the {\it cohomological support loci\/} by

$$V^i(X,A,\FF):=\{ P \in \Pic ^0(A) | h^i(X, \FF \otimes \pi^*P) > 0
\}.$$
In particular,
if $\pi={ {\alb }}_X \colon X \to \Alb (X)$,
then we simply write
$$V^i(X,\FF):=\{ P \in \Pic ^0(X) | h^i(X, \FF \otimes P) > 0 \}.$$

We sometimes write $V^i(\FF)$ instead of  $V^i(X,\FF)$ 
if no confusion is likely to arise.
We say that $P$ is a {\em general point}   of $V^i(X, A, \FF)$ 
if $V^i(X, A, \FF)$ is smooth
at $P$ and  $h^i(\FF\ot P)$ is  equal to the minimum 
of the function  $h^i(\FF\ot -)$
on  the component of
$V^i(X, A, \FF)$ that contains $P$.

In \cite{EL1}, Ein and Lazarsfeld illustrate various examples
in which the geometry of $X$ can
be recovered from information on the loci $V^i( \ox)$.
The geometry of $V^i(\ox )$
is governed by the following:

\begin{thm}\label{genvanish}{(Generic Vanishing Theorem)}

\noindent Let $X$ be a smooth projective
variety.  Then:
\begin{itemize}
\item[a)]$V^i(\ox )$ has  codimension
$\ge i-(\dim (X)-\dim( \Alb(X))$;

\item[b)] every irreducible component of $V^i(X,\ox)$ is
a translate of a sub-torus of\, $\PX$ by a torsion point;

\item[c)] let   $T$ be  an irreducible component of $V^i(\ox)$,  let
$P \in T$ be a point such that $V^i(\omega _X)$ is smooth at $P$, and let
$v
\in H^1(X, \OO
_X)\cong T_{P }\Pic ^0(X)$. If $v$ is not tangent to $T$, then the
sequence
$$H^{i-1}(X,\ox \otimes P) \stackrel {\cup v}  {\longrightarrow}
H^{i}(X,\ox \otimes P)  \stackrel {\cup v}  {\longrightarrow }
H^{i+1}(X, \ox \otimes P) $$
is exact. Moreover, if $P$ is a general point of $T$
and $v$ is
tangent to $T$ then
both maps   vanish;

\item[d)] if $X$ has maximal Albanese dimension, then there are inclusions:
$$ V^0(\ox )\supseteq V^1(\ox)
\supseteq \dots \supseteq V^n(\ox )=\{ \OO _X\}.$$

\end{itemize}
\end{thm}
\begin{proof}
Statement a) is Theorem 1 of \cite{GL1}. For statement b), 
the fact that the components of
$V^i(\ox)$ are translates of  abelian subvarieties follows 
from Theorem 0.1 of \cite{GL2} and
the fact  that the translation is by a torsion point follows 
from Theorem 4.2 of
\cite{Si}. Statement c) follows from Theorem 1.2 of \cite{EL1}. 
Statement d) is Lemma
1.8 of \cite{EL1}.
\end{proof}

\begin{rem}\label{generalizza}If $\pi\colon X\to A$ is a morphism to an
abelian variety, then
the loci
$V^i(X, A,\ox)$ satisfy properties analogous to a)\dots d) of Theorem
\ref{genvanish} (cf.
\cite{EL1}, Remark 1.6).
\end{rem}

The  results that follow are refinements of Theorem
\ref{genvanish} for
the case of a variety of maximal Albanese dimension with $\chi(\ox)=0$ or $1$.

\begin{prop}\label{V1isolated}
Let $X$ be a smooth projective variety of maximal Albanese dimension
with $\chi (\omega _X)=1$ and $q(X)\ge\dim (X)+1$.
If $P\in \Pic^0(X)$ is
an isolated point of $V^1(\ox)$, then $P=\OO_X$.
\end{prop}
\begin{proof}We write $n:=\dim (X)$
and $q:=q(X)\ge n+1$.  For  $n=1$ the claim follows by Serre duality, hence we
may assume $2\le n$. We will assume that $P\ne \OO _X$ and show that
$h^{n-1}(P^{-1}) =h^1(\omega _X \ot P)=0$. We remark that since $P\ne \OO _X$,
one has $h^0(P^{-1})=0$.
  Let
$k$ be the greatest integer such that
$h^k(\ox\ot P)>0$. Then by Theorem \ref{genvanish}, d),
$P$ is an isolated point of
$V^j(\ox)$ for $1\le j\le k$.
By Theorem \ref{genvanish}, c) and Serre duality, for every $0\ne v\in
H^1(\OO_X)$ the  complex:
$$0\to  H^0(P\inv)\overset{\cup v}{\to  }H^1(P\inv)\overset
{\cup v}{\to}\dots \overset
{\cup v}{\to}H^n(P\inv)$$
 is exact except at the last term.
As $v$ varies,  the above complexes fit together to give a complex of
vector bundles on
$\PPP := \PPP (H^1(\OO_X)\dual ))$ (see
\cite{EL1}, proof of Thm. 3, for a similar argument):
$$0\to H^0(P\inv )\ot \OO _{\PPP }(-n)\to H^1(P
\inv )\ot\OO _{\PPP }(-n+1)
\to  \dots\hskip2cm$$
$$\hskip4.5cm
\dots\to H^{n}(P\inv )\ot \OO _{\PPP }$$
which is again exact except at the last term, since exactness
can be checked fibrewise.
For $\OO _X\ne Q$ in an appropriate neighborhood of $\OO _X$,
one has that $h^i(P^{-1}\ot Q)=0$ for $i<n$,  hence $h^n(P^{-1}\ot Q)=\chi
(\omega _X\ot P \ot Q^{-1})=\chi (\omega _X)=1$.
It follows by \cite{GL2} Corollary 3.3,
that for every $0\ne v\in H^1(X,\OO _X)$, the
cokernel of the map $H^{n-1}(P^{-1})\overset{\cup v}{\to  }H^n(P^{-1})$
is $1-$dimensional. Therefore the cokernel of the map of vector bundles
$$H^{n-1}(P^{-1})\ot \OO _{\PPP }(-1)\lra
H^n(P^{-1})\ot \OO _{\PPP }$$
is a line bundle  $\OO _{\PPP }(d)$ with $d\geq 0$.
Let $\KK ^\bullet$ be the complex of vector bundles given by the
exact sequence
$$0\to H^0(P\inv )\ot \OO _{\PPP }(-n)\to H^1(P
\inv )\ot\OO _{\PPP }(-n+1)
\to  \dots\hskip2cm$$
$$\hskip1cm
\dots\to H^{n}(P\inv )\ot \OO _{\PPP }\to \OO _\PPP (d)\to 0.$$
For any $i$, one can  consider  the two
spectral sequences associated with the hypercohomology of the complex
$\KK ^\bullet \otimes \OO _\PPP (i)$.
Since the complex is exact, we have $'E^{s,t}_1=0$
for all
$s,t$, and thus
$'\!E^{s,t}_{\infty}=0$.  On the other hand, for $i=-1$,
one has that
$^{\prime\prime}\!E^{s,t}_1=0$ except for $s=0$ and $t=n+1$, in which case
$^{\prime \prime}\!E^{0,n+1}_1= H^0(\OO _\PPP (d-1))$ (recall that
$h^0(P\inv)=0$ by assumption). Comparing the two spectral sequences, one sees
that
$H^0(\OO _\PPP (d-1))=0$ and hence $d=0$. 
For $i=0$ the only non zero terms are\
$^{\prime \prime}\!E^{0,n}_1= H^n(P^{-1})\ot H^0(\OO _\PPP )$ and
$^{\prime \prime}\!E^{0,n+1}_1= H^0(\OO _\PPP )$, therefore
$h^n(P^{-1})=1$. In particular, the map
$H^n(P^{-1})\ot \OO _\PPP \to \OO _\PPP$ is an isomorphism.
For $i=1$ the only non zero terms are
$^{\prime \prime}\!E^{0,n-1}_1= H^{n-1}(P^{-1})\ot H^0(\OO _\PPP )$ ,
$^{\prime \prime}\!E^{0,n}_1= H^n(P^{-1})\ot H^0(\OO _\PPP (1))$ and
$^{\prime \prime}\!E^{0,n+1}_1= H^0(\OO _\PPP (1))$. The differential
$d_1\colon ^{\prime \prime}\!E^{0,n}_1\to\, ^{\prime \prime}\!E^{0,n+1}_1$
is an isomorphism and hence $^{\prime \prime}\!E^{0,n-1}_1= 0$ i.e.
$H^{n-1}(P^{-1})=0$.

\end{proof}

\begin{prop}\label{gencoho}
Let $X$ be a smooth projective  variety of maximal Albanese
dimension such that
$\chi(\ox)=0$. Let $T$ be a component of
$V^0(\ox)$ and let $W\subset H^1(\OO_X)$ be a linear subspace complementary to
the tangent space to $T$ (recall that $T$ is a translate of an abelian
subvariety of $\Pic^0(X)$). 
Let  $P\in T$ be a point such that $h^i(\ox\ot P)=0$ for $i>\dim (W)$ and such
that $P$ is a smooth point of $V^0(\ox)$. Then the natural
map induced by cup product:
$$H^0(\omega _X\ot P)\ot
\wedge ^{j} W\to H^j(\omega _X\ot P)$$
 is an isomorphism for every $j\ge 0$.
\end{prop}
Before proving the Proposition, we wish to point out that 
the assumptions on $P$ are
satisfied if $P\in T$ is general (cf. \cite{EL1}, proof of Theorem 3).
\begin{proof}Notice that by Theorem
\ref{genvanish},  a), the assumptions on
$X$ imply that
$V^0(\ox)$ is a proper subset of $\Pic^0(X)$.
Write $n:=\dim (X)$ and $w:=\dim (W)$. 
As explained in the proof of Theorem 3 of
\cite{EL1}, 
$T$ is a component of $V^w(\ox)$ and $T\not\subset V^i(\ox)$ for $i>w$.
By   Theorem  \ref{genvanish}, c), the
assumptions on $P$ imply that the complex
$D(v)$:
$$ \dots  H^{j-1}(\omega _X\ot P)
\overset {\cup v}  {\longrightarrow}
H^{j}( \omega _X\ot P)  \overset {\cup v}  {\longrightarrow }
H^{j+1}( \omega _X\ot P) \dots $$
is exact for all $v\in W\subset H^1(X,\OO _X)$.
As we have already seen in the proof of Proposition
\ref{V1isolated}, the complexes
$D(v)$ fit together to give an exact complex
$\KK ^{\bullet }$ of vector bundles on $\PPP:=\PPP(W)$:
$$0\to H^{0} ( \omega _X\ot P)\ot \OO _{\PPP} (-w)\to
H^1 (\omega _X\ot P)\ot \OO _{\PPP} (-w+1)\to\dots $$
$$\hskip3cm
\dots\to H^{n} ( \omega _X\ot P)\ot \OO _{\PPP}(-w+n)\to 0.$$
Similarly, there is an exact sequence $\KK ^{\bullet }_0$
of vector bundles on $\PPP$:
$$\dots\to
H^0 (\omega _X\ot P)\ot\wedge^{j-1} W\ot \OO _{\PPP}
(-w+j-1)\to H^0 (\omega _X\ot P)\ot\wedge^{j} W\ot \OO _{\PPP}
(-w+j)\to$$
$$\hskip3cm\to H^0 ( \omega _X\ot P)\ot \wedge ^{j+1}
W\ot \OO _{\PPP}(-w+j+1)\to\dots $$
There is a map of complexes $\KK ^{\bullet } _0\to \KK ^{\bullet }$
induced by
$$H^0 ( \omega _X\ot P)\ot \wedge ^{j}
W\ot \OO _{\PPP }(-w+j) \hookrightarrow
 H^0 ( \omega _X\ot P)\ot \wedge ^{j} H^1(X,\OO _X)
\ot \OO _{\PPP }(-w+j)
$$
$$\hskip3cm \to H^j ( \omega _X\ot P)\ot \OO _{\PPP }(-w+j)$$
Clearly for $j=0$ there is an isomorphism
$$H^0 ( \omega _X\ot P)\cong
 H^0 ( \omega _X\ot P)\ot\wedge ^0 W.$$
Proceeding by induction, assume $1\leq j\leq \dim (X)$ and
$$H^l ( \omega _X\ot P)\cong
H^0 ( \omega _X\ot P)\ot\wedge ^{l} W$$ for all $l<j$.
Tensoring by $\OO _{\PPP}(-j)$ and
taking cohomology, one gets:
\smallskip
$$
\CD
\dots & & & & \dots\\
@VVV & & @VVV & \\
H^0\hskip-.02cm
(\omega  _X \ot P)\hskip-.02cm\ot \hskip-.08cm\wedge ^{j
 -1}W
\hskip-.02cm\ot  \hskip-.02cm H^{w - 1}\hskip-.02cm(\OO _{\PPP }
(-w-1))\hskip-.12cm& @>{\sim}>>&\hskip-.12cm
H^{j-1}\hskip-.02cm(\omega _X \ot P)
\hskip-.02cm\ot \hskip-.02cm H^{w-1}\hskip-.02cm(\OO _{\PPP }(-w\hskip-.02cm
-\hskip-.02cm 1))\\
@VVV & & @VVV & \\
H^{0}\hskip-.02cm(\omega  _X\ot P)\ot \wedge ^{j} W
\ot  H^{w-1}\hskip-.02cm(\OO _{\PPP }(-w))&\hskip-.12cm
@>>>\hskip-.12cm&H^{j}\hskip-.02cm
(\omega _X\ot P) \ot  H^{w-1}\hskip-.02cm(\OO _{\PPP }(-w)) \\
@VVV & & @VVV & \\
0 & & & & 0
\endCD
$$
If the vertical rows are exact, then
the required isomorphism follows from the five lemma .
 Consider the spectral sequence associated to the complex
$\KK ^{\bullet} \ot \OO _{\PPP }(-j)$. We have $'E_1^{s,t}=0$ for any
$s,t$, thus $'E_{\infty}^{s,t}=0$.
On the other hand for $1\leq j \leq \dim (X)$ one has
$''E_1^{s,t}=0$ if $s\ne w-1$ and $s\ne 0$. If $s=0$ then $''E_1^{0,t}=
H^t(\omega _X \ot P)\ot H^0(\OO _{\PPP }(-w+t-j))$.
For $0\leq t \leq w$ this group is $0$ as $-w+t-j<0$, 
while for  $t>w$ it is $0$ since $h^t(\omega _X \ot P)=0$ by assumption.
Comparing spectral sequences one sees that the sequence
$$...\longrightarrow
H^{j-1}(\omega _X \ot P)\ot H^{w-1}(\OO _{\PPP }(-w-1))\longrightarrow
H^{j}(\omega _X \ot P)\ot H^{w-1}(\OO _{\PPP }(-w))\longrightarrow
...$$
is exact. The exactness of the first vertical row can be shown in the same way.
\end{proof}
The following Corollary is due to Ein--Lazarsfeld (cf. \cite{EL2}).
\begin{cor}\label{Vzero}Let $X$ be a smooth projective  
variety of maximal Albanese
dimension such that
$\chi(\ox)=0$. If $P\in \Pic^0(X)$ is an isolated point 
of $V^0(\ox)$ then $P=\OO_X$.
\end{cor}
\begin{proof}Write as usual $n:=\dim X$, $q:=q(X)$. 
By Proposition \ref{gencoho} we have
$H^q(\ox\ot P)\cong
\wedge^q H^1(\OO_X)\ot H^0(\ox\ot P)\ne 0$. This implies $q\le n$, hence $q=n$,
since $X$ has maximal Albanese dimension. 
In particular we have $0\ne h^n(\ox\ot
P)=h^0(P\inv)$, implying $P=\OO_X$.
\end{proof}

 In  some cases it is possible to obtain more
precise information on the linear systems
$|mK_X+P|$ for
$m>1$ and $P\in V^0(\omega _X^m).$

\begin{prop}\label{Pm}Let $X$ be a smooth projective variety and
let $f\colon X\to Y$ be the Iitaka fibration (assume $Y$ smooth). Fix
$m\ge 2$  and
$Q\in
\Pic^{\tau}(X)$.
Then:

$$h^0(mK_X+ Q+ P)=h^0(mK_X+Q)\,\,\, \text{for all}\,\, P\in \Pic ^0(Y).$$
\end{prop}
\begin{proof}
The statement is trivial if $\kappa(X)=\dim (Y)=0$,
hence we may assume $\kappa(X)>0$.

\noindent{\bf Step 1.}
{\em  For all $P\in \Pic ^0(Y)$,  $\kappa
(K_X+P)\geq \kappa(X)$.}

 Arguing as in the proof of
\cite[Lemma 2.1]{CH1}, one can show that    if
$m\ge 2$ then $h^0(mK_X)=h^0(mK_X+P)$ for all  $P\in \Pic^{\tau}(Y)$. 
Since $\Pic ^\tau (Y)$
is dense in $\Pic^0(Y)$, by semicontinuity one has that 
$h^0(m(K_X+P))=h^0(mK_X+mP)\geq
h^0(mK_X)$ for every $P\in
\Pic^0(Y)$.
\medskip

\noindent{\bf Step 2.} {\em Let
$P\in \Pic ^0(Y)$. We have    $h^0(mK_X+
P+Q)=h^0(mK_X+P+R+ Q)$ for  all $R\in \Pic^{\tau}(Y)$}

This follows by a procedure analogous to \cite[Lemma 2.1]{CH1}.
We illustrate this for $m=2$. Write $S:=\Alb(Y)$ and denote by $\pi\colon X\to S$ the
composition of $f$ and of the Albanese map of $Y$.  Fix $H$ an ample line bundle on
$S$ and $R\in \Pic^{\tau}(Y)$. By Step 1,   for $r$ sufficiently big and divisible
 we may pick a divisor:
$$B\in |r(K_X+P)-\pi^*H|= |r(K_X+P+Q+R)-\pi^*H|=|r(K_X+P+Q)-\pi^*H|.$$
Possibly replacing $X$ by an appropriate birational model, we may assume  
that $B$
has normal crossings support and that (as in \cite[Lemma 2.1] {CH1})
$$\lfloor \frac {B}{r}\rfloor \subset Bs |2K_X+P+Q| \cap Bs |2K_X +P+Q+R|.$$
Let $$L:=K_X-\lfloor \frac {B}{r}\rfloor \equiv \frac {\pi^*H}{r}+\{ \frac {B}
{r}\}$$
i.e. $L$ is numerically equivalent to the pull back of a nef and big 
$\QQ-$divisor plus a klt
divisor. Then comparing base loci as in \cite[Lemma 2.1]{CH1}
$$h^0(2K_X+P+Q)=h^0(K_X+L+P+Q)=$$
$$=h^0(\pi_*(\omega _X\ot L\ot P\ot Q))=
\chi (\pi_*(\omega _X\ot L\ot P\ot Q))=$$
$$=\chi (\pi_*(\omega _X\ot L\ot P\ot Q\ot R))=
h^0(\pi_*(\omega _X\ot L\ot P\ot Q\ot R))=$$
$$ =h^0(K_X+L+P+Q+R)=h^0(2K_X+P+Q+R),$$
where the third and the fifth equality  follow from Theorem \ref{omegavan},
b).
\medskip

\noindent{\bf Step 3.} {\em For all
$P\in \Pic ^0(Y)$, $h^0(mK_X+ P+Q)=
h^0(mK_X+Q)$}.

Let $M$ be the maximum of the function $h^0(mK_X+Q+P)$ for  
$P\in \Pic^0(Y)$ and let $P_0\in
\Pic^0(Y)$ be such that $h^0(mK_X+Q+P_0)=M$. By Step  2, 
$h^0(mK_X+ Q+P_0+R)=M$ for all
$R\in \Pic^{\tau}(Y)$. Since $\Pic^{\tau}(Y)$ is dense in
$\Pic^0(Y)$,  by semicontinuity we have $h^0(mK_X+Q+P)\ge M$ 
for all $P\in \Pic^0(Y)$, hence
$h^0(mK_X+Q+P)=M$ for all $P\in \Pic^0(Y)$.
\end{proof}
In \S \ref{P3=2} we will need following result, which   
is due to Ein and Lazarsfeld.
\begin{lem}\label{lemmaD}Let $X$ be a variety such that $\chi(\ox)=0$
and such that
$\alb_X\colon X\to
\Alb(X)$ is surjective and  generically finite.    Let
$T$ be an irreducible component of
$V^0(\omega _X)$,  and let
$\pi _E\colon X\to E:=\Pic^0(T)$ be the morphism induced by the
map
$\Alb(X)=\Pic^0(\Pic^0(X))\to E$ corresponding  to the inclusion
$T\hookrightarrow\Pic^0(X)$.

Then there exists a divisor $D\prec R:=Ram (\alb _X) =K_X$, vertical
with respect to $g:=\pi_E\circ \alb _X$, 
such that for general $P\in T$, $F:=R-D$ is a fixed
divisor
of each of the linear series $|K_X+P|$.
\end{lem}
\begin{proof} This is \cite{EL2}, Lemma D. We thank Ein and Lazarsfeld
for allowing us to reproduce their proof.

Let $n=\dim (X)=q(X)$, and let $k$ be the codimension of $T$ in $\Pic ^0(X)$.
Choose a basis
$$\omega _1,..., \omega _n \in H^0(X,\Omega ^1_X)$$
such that $\omega _{k+1},...,\omega _n$ are pull backs
under $g$ of a basis for $H^0(E,\Omega ^1_E)$. In particular we may assume that
$\omega _{k+1},...,\omega _n$ are conjugate to a basis of the tangent space
to $T$.
There is a homomorphism:
$$\alpha: \OO _X =g^* \Omega ^{n-k}_E\lra \Omega ^{n-k}_X$$
defined by the section $\omega _{k+1}\wedge ... \wedge \omega _n
\in H^0(X,\Omega ^{n-k}_X)$. $\alpha$ vanishes at a point $x\in X$
if and only if $g$ is not smooth at $x$. Let $D\subset X$ be the divisor 
along which $\alpha$ vanishes. There is a diagram:

$$
\CD
 &@. &\OO _X&@.& &@.& & @. &\\
@. & & @VV{D}V    & &  @.& & @.& & @. \\
0 & @>>>&\OO _X(D) &@>{\bar{\alpha }}>>&\Omega ^{n-k}_X & @>>> 
& coker(\bar{\alpha }) & @>>> &0.\\
\endCD
$$
An easy local computation shows that $coker(\bar{\alpha })$ is torsion
free. Fix a general point $P\in T$.
By the proof of Theorem 3 of \cite {EL1} (see also Proposition \ref{gencoho})
one sees that $H^0(X,\Omega _X ^{n-k}\ot P) \ne 0$.
We claim that $\bar {\alpha}$ induces an isomorphism:
$$ H^0(X,\OO _X(D)\ot P)=H^0(X,\Omega _X ^{n-k}\ot P).$$
Fix a non-zero section $\eta \in H^0(X,\Omega _X ^{n-k}\ot P)$. Since
$coker(\bar{\alpha })$ is a torsion free sheaf, it suffices to show
that at a general poit $x\in X$, $\eta (x)$ lies in the subspace of
$\Omega ^{n-k}_X\ot P(x)\cong \Omega _X^{n-k}(x)$
spanned by $\omega _{k+1}(x)\wedge ... \wedge \omega _n(x)$. By Theorem
\ref{genvanish}
for all $k+1\leq i \leq n $, one has 
$\eta \wedge \omega _i=0\in H^0(X,\Omega _X^n\ot P)$. Since $\omega _1(x),...,
\omega _n(x)$ span $T_x^*(X)$ at a general point, $\eta (x)$ must be a multiple
of $\omega _{k+1}\wedge ... \wedge \omega _n$.

Consider finally the commutative diagram:
$$
\CD
P &@>{R}>> &\Omega ^n_X \ot P\\
@V{D}VV  & & @AA{\wedge \omega _1\wedge ... \wedge \omega _k}A \\
\OO _X(D)\ot P & @>{\bar {\alpha }}>>&\Omega ^{n-k}_X \ot P. 
\endCD
$$
The top row is multiplication by the ramification divisor $R=\{
\omega _1\wedge ... \wedge \omega _n=0\}$, therefore $D\prec R$.
As we have seen above, $\bar{\alpha }$ gives an isomorphism on
global sections. By the proof of Theorem 3 of \cite{EL1}
(cf. Proposition \ref{gencoho}) one sees that the right 
hand side vertical homomorphism
also induces an isomorphism on
global sections. Therefore
$$H^0( \OO _X(D)\ot P )=H^0(K_X\ot P).$$
In other words $F=R-D$ is a fixed divisor of the linear series
$|K_X+P|$.
\end{proof}

Finally, we often use the following observation to give lower bounds
on the dimension of the
series $|mK_X+P|$.
\begin{lem}\label{claimA}Let $X$ be a smooth projective variety,
let $L$ and $M$ be line bundles on $X$, and let  $T\subset \Pic ^0(X)$
be an irreducible subvariety of dimension $t$. If for all $P\in T$,
$\dim |L+ P|\geq
a$ and $\dim |M-P|\geq b$, then $\dim |L+ M|\geq a+b+t$
\end{lem}
\begin{proof}
Denote by ${\mathcal P}$ the restriction to $X\times T$ of the Poincar\'e line
bundle on $X\times \Pic^0(X)$ and by $p_i$, $i=1,2$, the projections of
$X\times T$
on the
$i-$th factor. Set $\V_1:={p_2}_*({\mathcal P}\ot p_1^*M)$ and $\V_2:=
{p_2}_*({\mathcal P\inv}\ot L)$ and denote by $\rho_i$\  the generic rank
of $\V_i$,
$i=1,2$.
Denote by $\psi\colon \V_1\ot \V_2\to {p_2}_*(p_1^*(L\ot M))$
the sheaf map induced by
the multiplication map $$(p_1^*L\ot {\mathcal P})\ot
(p_1^*M\ot {\mathcal P\inv})\to
p_1^*(L\ot M)$$ on $X\times T$. The sheaf ${p_2}_*(p_1^*(L\ot M))$ is
isomorphic to  the
trivial bundle $\OO_T\times H^0(X, L\ot M)$.
In addition, there  exists a nonempty open set $T_0\subset T$ 
such that $\V_1$  and
$\V_2$ are locally free on $T_0$.   
We denote by $V_1$, $V_2$ the vector bundles
associated to the restriction of
$\V_1$, $\V_2$ to $T_0$. For each $P\in T_0$ there are natural
identifications
$V_{1, P}\cong H^0(X, L+P)$ and $V_{2, P}\cong H^0(X, M-P)$. In particular we
have $\rho_1\ge a+1$ and $\rho_2\ge b+1$. Composing the natural bilinear map
$V_1\times_{T_0} V_2\to V_1\ot V_2$ with the restriction of $\psi$, 
we obtain a morphism
$V_1\times_{T_0} V_2
\to T\times H^0(X,L\ot M)$. In turn, passing to the projectivized bundles and
composing with the projection $T_0\times |L\ot M|\to |L\ot M|$, this induces a
morphism
$\phi\colon \PPP(V_1)\times_{T_0}\PPP(V_2)\to |L\ot M|$ 
such that the restriction of
$\phi$ at the fibers over $P\in T_0$ corresponds to the natural map of linear
systems $|L+P|\times |M-P|\to |L\ot M|$. The map $\phi$ has finite fibers, 
since each  
element of $|L\ot M|$ can be decomposed as the sum of two effective divisors
in a finite number of ways. It follows that the dimension of 
$|L\ot M|$ is greater
than or equal to
$\dim\PPP(V_1)\times_{T_0}\PPP(V_2)=\rho_1+\rho_2+t-2\ge a+b+t$.
\end{proof}


\section{Birationality criteria}\label{birationality}
In this section we exploit the techniques of
\ref{mukai}, \ref{dualizing} and \ref{supportloci} to give
criteria for the birationality of maps between
varieties of maximal Albanese dimension.

\begin{thm}\label{bir1} Let $f\colon X\to Y$ be a generically finite
morphism of smooth  projective
 varieties of maximal Albanese dimension.
If the induced maps
$H^i(Y,\omega _Y\ot P)\to  H^i(X,f_*\omega _X\ot P)$
are isomorphisms for all $i\geq 0$ and all $P\in \Pic ^0(Y)$,
 then $f$ is birational.
\end{thm}
\begin{proof} Denote by $j\colon \omega _Y \to f_*(\omega _X)$ the natural
inclusion. Pushing forward to $\Alb(Y)$, one has an  
inclusion ${\alb _Y}_*(\omega
_Y)\hookrightarrow {\alb _Y}_*(f_*(\omega _X))$.
By Theorem \ref{omegavan}, c), the
sheaves $R^i{\alb _Y}_*(f_* \omega _X)$ and
$R^i{\alb _Y}_*\omega _Y$ vanish for all $i>0$, hence for all 
$P\in Pic^0(Y)$ the Leray
spectral sequence gives natural  isomorphisms 
$H^i(Y, \omega_Y\ot P)\cong H^i(A,{\alb
_Y}_*(\omega _Y)\ot P)$ and $H^i(Y,f_*\omega _X\ot P)$
$\cong H^i(A,{\alb _Y}_*(f_*(\omega _X))\ot P)$. Thus for all
$P\in
\Pic ^0(Y)$ and
$i\geq 0$ we have isomorphisms
$$H^i(A,{\alb _Y}_*(\omega _Y)\ot P)\stackrel{\sim}{\to}
H^i(A,{\alb _Y}_*(f_*(\omega _X))\ot P)\ \ \ .$$
By Corollary \ref{inclusion}, the sheaves ${\alb _Y}_*(\omega _Y)$ and
${\alb _Y}_*(f_*(\omega _X))$ are isomorphic and, in particular, they have
the same rank at the generic point of $\alb_Y(Y)$. Thus  the
degree of
$X\to  Y$ is equal to 1, i.e.
$f$ is birational.
\end{proof}

Let $f\colon X\to Y$ be a surjective
morphism of smooth  projective
 varieties of the same dimension.
Let $X\to V$ and $Y\to W$ be birational models of the
respective Iitaka fibrations.
We have an inclusion of
linear series
$$|mK_Y|\stackrel{f^*}{\lra} |mf^*K_Y|\stackrel{+mK_{X/Y}}{\lra} |mK_X|.$$
Since $Y\to  W$ is defined by sections of $|mK_Y|$ for $m$ sufficiently
big and divisible, and since the pull-backs of these sections
correspond
sections of $|mK_X|$, one sees that the induced morphism $X\to  W$
factors through a rational map $V\to  W$.
If $P_m(X)=P_m(Y)$ for $m$ sufficiently
big and divisible, then $V\to W$ is birational.
In particular, if  $Y$ is of general type, then $X$ is birational to $Y$.

It would be interesting to obtain effective versions of this result.
In general this seems to be a very complicated problem
and it is not clear
what bounds to expect. However when $Y$ is of maximal Albanese dimension
we obtain some satisfactory results.

\begin{thm} Let $f\colon X\to Y$ be a dominant
morphism of smooth $n$-dimensional  projective
 varieties, $Y$ of maximal Albanese dimension, and let $X\to V$, 
$Y\to W$ be the
Iitaka fibrations of $X$, respectively $Y$. If $P_2(X)$ $=P_2(Y)$, 
then the induced map
$V\to  W$ has connected fibers.  In particular,
if $Y$ is of general type, then
$f$ is birational.
\end{thm}

\begin{proof}
If $\kappa(Y)=0$, then $W$ is a point and the claim is of course true. 
Therefore we
may assume $\kappa(Y)>0$.

By \cite {CH3}, \S 1.1, we may assume that $V$ and $W$ are smooth and that  
the Iitaka fibrations of $X$ and
$Y$ fit in a commutative diagram where $A:=\Alb (Y)$ and $S:=\Alb (W)$
$$
\CD
X &@>f>>&Y&@>a>>&A \\
@VVV & & @VV{p_Y}V    & &  @VVV \\
V & @>>>&W&@>s>>&S.\\
\endCD
$$

Consider the induced maps $p_X\colon X\to W$, $\pi _X \colon
X\to S$ and
$\pi _Y
\colon Y\to S$. The claim will follow if we show that $p _X$
has connected fibers. Fix $H$
ample on
$S$. The map $s\colon W\to  S$ is generically finite onto its
image by Proposition
\ref{albanese}, a), hence
$s^*H$ is nef and big on $W$. For  $m\gg 0$,
the linear system $|mK_Y-\pi _Y^*H|$ is
nonempty. Let
$B$ be a general member of
$|mK_Y-\pi _Y^*H|$ and let $\tilde{B}:=f^*B+mK_{X/Y}$.
We may assume that both  $B$ and
$\tilde{B}$ have normal crossings support.
Define
$$L_Y:= \omega _Y (-\lfloor B/m\rfloor )\equiv {\pi }_Y^*(H/m)+\{ B/m\}$$
 and
$$L_X:=\omega _X (-\lfloor \tilde{B}/m\rfloor )=f^* \omega _{Y}(-\lfloor (f^*B)
/m\rfloor )\equiv {\pi} _X^*(H/m)+\{ \tilde{B}/m\}.$$

\noindent{\bf Step 1.} {\em  $K_{X/Y}-\lfloor (f^*B)/m\rfloor +f^*\lfloor
B/m\rfloor$ is effective.}
\smallskip

Let $B/m=\sum_1^s b_iB_i$, where the   $B_i$ are distinct prime divisors.
If the $b_i$ are integers, then $\lfloor (f^*B)/m\rfloor
=f^*\lfloor B/m\rfloor $.
So it is enough to consider the case $0\le b_i<1$, i.e. 
$\lfloor B/m\rfloor=0$. Let
$P\in Y$ be a point such that $P\in B_i$ for $1\le i \le s$ and let   
$(y_1,...,y_n)$
be local coordinates centered 
in $P$ such that $y_i$ is a local equation for $B_i$, $i=1\dots s$. 
Let $Q$ be a
point such that $f(Q)=P$ and let $E$ be   a component of $\tilde{B}$
containing $Q$.  Choose
local coordinates $(x_1\dots x_n)$ around $Q$ such that $x_1$
is a local equation
for $E$. Assume that $E$ is a component of $f^*B_i$ for $i=1\dots t$,
so that for
$i=1\dots t$  we have $f^*y_i=x_1^{n_i} \epsi_i$, with $n_i>0$ and $\epsi_i$ a
regular function that does not vanish identically on $E$.
A local equation for $K_{X/Y}$ around $p$ is given by the determinant of the
Jacobian matrix $(\frac{\partial f^*y_i}{\partial x_j})$, which is easily
seen to
vanish  on $E$ to order at least $(\sum_1^tn_i)-1$. So the coefficient of
$E$ in
$K_{X/Y}-\lfloor (f^*B)/m\rfloor +f^*\lfloor B/m\rfloor $ is greater than
or equal to
$\left( \sum_1^t(n_i-\lfloor n_ib_i\rfloor)\right)-1\ge 0$.
\medskip

\noindent{\bf Step 2.} {\em There is a map of sheaves
$\omega _Y\ot L_Y\to f_*(
\omega _X\ot L_X)$
inducing an isomorphism on global sections.}
\smallskip

Since $Y$ is of Albanese general type, we have $P_m(Y)>0$ for every $m\ge 1$.
After replacing $X,Y$ by appropriate birational models,
we may assume that
$$|2K_Y|=F_2+|M_2|$$
where  $|M_2|$ is free and $B+F_2$ has normal crossings,  $B\in |mK_Y-\pi _Y^*
H|$ being  the divisor chosen before.

We wish to define a new divisor $B'\in |m'K_Y-\pi _Y^*H|$ such that
$\lfloor B'/m'\rfloor \prec F_2$. To this end, pick $D=F_2+D'\in |2K_Y|$
such that $D+B$ has normal crossings support and $D'$ is smooth, not contained
in the support of $B$. Let $m':=m+2l$ and $B':=B+lD\in |m'K_Y
-\pi ^*_Y H|$. We may write $D=\sum d_i B_i +D'$ and $B=\sum b_iB_i$.
For all $l\gg 0$ we have that the multiplicity of $B'$ along $B_i$
satisfies
$$mult _{B_i}\lfloor B'/m'\rfloor =\lfloor \frac {b_i +ld_i}{m+2l}\rfloor
\leq \lfloor \frac {b_i}{2l}+\frac {d_i}{2}\rfloor \leq d_i
=mult _{B_i}F_2.$$
We now replace $B$ by $B'$. It follows that $h^0(\omega _Y \ot L_Y)=P_2(Y)$.

By Step 1, there is an injection of sheaves $f^*(\omega _Y \ot L_Y)
\to \omega _X \ot L_X$,
corresponding to the effective divisor $K_{X/Y}-\lfloor (f^* B)/
m\rfloor +f^*\lfloor B/m\rfloor$. Pushing forward to $Y$, one gets an 
inclusion $\omega_Y\ot
L_Y\to f_*(\ox\ot L_X)$. Since $h^0(\omega _Y\ot L_Y)=P_2(Y)=
P_2(X)\geq h^0(f_*(K_X+L_X))$,
the corresponding map  on global sections is an isomorphism.
\medskip

\noindent{\bf Step 3.} $(\pi _X)_*(\omega _X \ot L_X)=
(\pi _Y)_*(\omega _Y \ot L_Y).$

By step 2, there is an injection of sheaves $\omega _Y \ot L_Y
\to f_*(\omega _X \ot L_X)$ inducing an isomorphism on global sections.
Pushing forward via $\pi _Y$ we obtain an exact sequence
of sheaves on $S$
$$0\to {\pi _Y }_*(\omega _Y \ot
L_Y)\to{\pi _X }_*(\omega _X \ot L_X) \to \QQQ \to 0.$$
By Theorem \ref{omegavan}, b), for all $i>0$, $P\in \Pic ^0(S)$,
$$H^i(S,{\pi _Y}_* (\omega _Y \ot L_Y)\ot P)=H^i(S,{\pi _X}_* (\omega _X \ot
L_X)\ot P)=0.$$ Therefore $h^i(\QQQ \ot P)=0$ for all $i>0$ and
$P\in \Pic ^0(S)$
and we have $h^0(\QQQ\ot P)=\chi(\QQQ\ot P)=\chi(\QQQ)$.  By Step 2,
the map $$H^0({\pi _Y }_*(\omega _Y \ot
L_Y))\to H^0({\pi _X }_*(\omega _X \ot L_X))$$ 
is an isomorphism.  Therefore we
have
$0=h^0(\QQQ )=h^0(\QQQ \ot P)$  and hence $\QQQ =0$ by Lemma \ref{nonvanish}.
\medskip

\noindent{\bf Step 4.} {\em $F_{X/W}$ is connected.}

Since $s\colon W\to S$ is generically finite onto its image,
it follows from Step 3
that the generic rank of ${{p} _X}_*(\omega _X \ot L_X)$
is equal to the generic rank of ${{p} _Y}_*(\omega _Y \ot L_Y)$.
Let $F_{Y/W}$ be a general
fiber of $Y\to W$. Recalling that $\omega_Y\ot L_Y$ is effective and that
$\kappa(F_{Y/W})=0$, we have
$0<h^0(F_{Y/W},
\omega _Y
\ot L_Y|_{F_{Y/W}})\le P_2(F_{Y/W})\le 1$. This shows that the generic
rank of ${{p}
_Y}_*(\omega _Y \ot L_Y)$ is equal to
$1$. Since
$\omega _X
\ot L_X$ is effective, the generic rank
of ${{p} _X}_*(\omega _X \ot L_X)$ is greater than  or equal to the
number of connected components of the generic geometric fiber.
Therefore, $F_{X/W}$ is connected.
\end{proof}
\section{Abelian varieties and theta divisors}
In this section  we
apply the results of \S \ref{birationality} to give some new
characterizations of abelian varieties and theta divisors.
We will need the following result:
\begin{prop}\label{iso}
Let $f\colon X\to A$ be a morphism from a smooth projective 
variety of dimension $n$ to an
abelian variety of dimension $n+1$. If $f^*\colon 
H^0(A, \Omega_A^n)\to H^0(X,\Omega^n_X)$
is an isomorphism, then  
$f^*\colon H^0(A,
\Omega_A^i)\to H^0({X},\Omega^i_{{X}})$ 
is an isomorphism for $0\le i\le n$.
\end{prop}
\begin{proof} 
The assumption 
implies that $\bar{X}:=f(X)$ is a divisor of $A$, hence $f$ is
generically finite onto its image. In addition, $\bar{X}$ is of general type, 
since
otherwise, by \cite[10.9]{Ue}, $\bar{X}$ would be a  
pull-back of a divisor from a quotient
$\bar{A}$ of $A$ of dimension $d\le n$. In that case, 
if $\omega_1, \dots ,\omega_d$ are the
pull-backs  of the elements of  a basis of 
$H^0(\bar{A},\Omega^1_{\bar{A}})$, then the
restriction of
$\omega_1\wedge\dots\wedge
\omega_d$\  to $\bar{X}$ is zero, contradicting the assumptions. 
Notice that the map
$f^*\colon H^0(A, \Omega_A^i)\to H^0(X,\Omega_X^i)$ is injective 
for every $i\le n$ by
the assumption that it is injective for $i=n$. 
The proof of Corollary 3.11 of \cite{Mo} gives
$h^0(\Omega_A^i)=h^0(\Omega^i_{ {X}})$ and thus 
$H^0(A, \Omega_A^i)\cong H^0({{X}},\Omega_{{X}}^i)$
for $0\le i\le n$.
\end{proof}

 The  result that follows
generalizes Theorem \ref{theta}.
\begin{prop}\label{thetacohomology}
Let $f\colon X\to A$ be a morphism from a smooth variety
of dimension $n$ to an abelian variety of dimension $n+1$. If
$h^0(X,\omega _X \ot f ^*P)=1$ for all $\OO _A \ne P\in \Pic ^0(A)$
and the map $$f ^*\colon H^0(A,\Omega ^{n}_A)\to H^0(X,\Omega ^{n}_X)$$ is an
isomorphism, then $f$ is birational onto its image $\bar {X}$ and
$\bar {X}$ is a principal polarization.
\end{prop}
\begin{proof}

Let $\tilde {X}$ be a desingularization of $\bar {X}$.  We
may  assume that $f$ factors through $\tilde {f}\colon X\to\tilde {X}$ and
$\phi \colon
\tilde {X}\to A$.

\noindent{\bf Step 1.} {\em  $\bar{X}$ is a principal polarization.}

Notice first of all that,  by the
injectivity of
$$f ^*\colon H^0(A,\Omega ^{n}_A)\to H^0(X,\Omega ^{n}_X),$$ $\bar{X}$ is a
divisor and  generates $A$. In particular, $X$ has maximal Albanese dimension.
By Theorem \ref{theta}, to prove that
$\bar{X}$ is a principal polarization it is  enough to show that:

  a)
$h^i(\omega _X)=\binom {n+1}{i+1}$ for all
$0\leq i\leq n$;   

 b) $h^i(\omega _X \ot f ^* P)=0$ for all
$i\ge 1$ and $\OO_A \ne P\in \Pic^0 (A)$.

\noindent
Condition a) follows directly from Proposition \ref{iso}.

To prove condition b), assume by way of contradiction that $V:=\cup_{i\ge 1}
V^i(X, A,\omega_X)$ contains a  point $P\ne \OO_X$. 
By Theorem \ref{genvanish}, a), and
Remark \ref{generalizza}, $V$ is a proper subset of $A$. Let
$T$ be a component of $V$, let  $P\in T$ be a general point and let  
$v\in H^1(A,\OO
_A)$ be a vector that is not tangent to
$T$ at $P$. By  Theorem \ref{genvanish}), c) and Remark \ref{generalizza},
the complex $(H^i(X,\omega_X\ot P), \wedge v)$ is exact
for $i\geq 1$. In particular,  we have
$\chi ( \omega _X\ot P)\leq h^0(\omega _X\ot P)$,
with equality holding iff the map
$H^0(X,\omega  _X\ot P )\stackrel{\wedge v}{\to}
H^1(X,\omega  _X\ot P )$ is zero.
Since $X$ is of maximal Albanese dimension, it is easy to see
that there exists $\sigma \in H^0(A, \Omega ^1_A)$ such that the
map $H^0(X,\Omega ^{n-1} _X\ot P)\stackrel{\wedge \si}{\to}
H^0(X,\omega  _X\ot P )$ is non-zero. 
If we denote by $v\in H^1(A,\OO _A)$ the conjugate of $\si $, 
then by Theorem \ref{genvanish} c) 
 $v$ is not tangent to $V^1(\ox )$ at $P$.
By Hodge conjugation and Serre duality with 
twisted coefficients (cf. \cite{GL1}, 2.5), the
map
$H^{0}(X,
\omega  _X\ot P)\stackrel{\wedge v}{\to} H^{1}(X,
\omega  _X\ot P)$ is non-zero. Therefore,
$\chi ( \omega_X\ot P)<h^0(\omega _X\ot P)=1$. On the other hand,  $\chi (
\omega_X\ot P)=\chi(\omega_X)$, and $\chi(\omega_X)=1$ as observed above. 
Thus we have
a contradiction and the proof of Step 1 is complete.

\medskip

\noindent{\bf Step 2.}\ \ {\em The map $f$ is birational.}

Since $f$ factors through $\phi\colon \tilde{X}\to A$,
 it is clear that $\phi^*\colon H^0(A,\Omega^{n} _A)\to H^0(\tilde {X}, 
\Omega^{n}
_{\tilde {X}})$ is an isomorphism. Then Proposition \ref{iso} implies  that
the maps $\phi ^*\colon H^0(A,\Omega ^{i}_A)
\to H^0(\tilde {X},\Omega ^{i}_{\tilde {X}})$ are
isomorphisms for $0\leq i\leq n$. In particular one has
$\chi(\omega_{\tilde{X}})=1$. In addition, for every $P$ the map
$\tilde{f}^*\colon H^0(\tilde{X}, \omega_{\tilde{X}}\ot P)\to H^0(X,\omega_X\ot
f^*P)$ is injective and we have $h^0(\tilde{X}, \omega_{\tilde{X}}\ot P)\le
h^0(X,\omega_X\ot f^*P)=1$. Now  the same argument as in   
Step 1 can be used to
show that
$h^i(\tilde{X},\omega_{\tilde{X}}\ot P)=0$ for $P\ne\OO_A$ and $i>0$. Thus for
every $P\ne \OO_A$ we have $h^0(\tilde{X},\omega_{\tilde{X}}\ot
P)=\chi(\omega_{\tilde{X}}\ot P)=\chi(\omega_{\tilde{X}})=1$.
It follows that for $P\ne \OO_A$ the map 
$H^0(\tilde{X},\omega_{\tilde{X}}\ot P)\to
H^0(\tilde{X},\tilde{f}_*\omega_X\ot P)=H^0(X,\omega_X\ot f^* P)$ is an
isomorphism, being a non trivial map of $1-$dimensional vector spaces. So the
assumption of Theorem \ref{bir1} is satisfied for $P\ne \OO_A$. For
$P=\OO_A$ it follows from the fact that
$H^0(\tilde{X},\Omega_{\tilde{X}}^i)\to H^0(X,
\Omega^i_X)$ is an isomorphism  by using Hodge conjugation and Serre
duality. Thus $f$ has degree 1 by Theorem \ref{bir1}.  \end{proof}

\begin{cor}\label{pgq} Let $X$ be a smooth projective variety such that
$$P_1(X)=q(X)=\dim (X)+1.$$
If $X$ is of maximal Albanese dimension and $\Alb (X)$ is simple, then
$X$ is birational to a theta divisor in $\Alb (X)$.
\end{cor}
\begin{proof}We write $n=\dim (X)$,
 and $A:=\Alb(X)$. Consider the map $H^0(A,
\Omega^n_A)\to H^0(X,\Omega^n_X)$. If it is not
an isomorphism, then there exist $1-$forms
$\omega_1,\dots
,\omega_n\in H^0(A,\Omega^1_A)=H^0(X,\Omega^1_X)$
such that $\omega_1\wedge\dots\wedge
\omega_n=0$ on $X$. By the Generalized Castelnuovo - de Franchis
Theorem (cf. for instance
\cite{fab}, Theorem 1.14) it follows that there exists a fibration
$p\colon X\to Y$ ($\dim (Y)<\dim (X)$) such that
$\omega_1\dots
\omega_n$ are pull-backs from $Y$.  Notice that $n\le q(Y)<q(X)$,
since $X$ has maximal
Albanese dimension. Then the induced morphism  $A \to \Alb(Y)$ is
surjective  and is not an isogeny, contradicting 
the assumption that $A$ be simple.
Thus $H^0(A, \Omega^n_A)\to H^0(X,\Omega^n_X)$ is an isomorphism and,  
by Proposition
\ref{iso},
$H^0(A, \Omega^i_A)\to H^0(X,\Omega^i_X)$ is also an isomorphism for all $i\leq
n$. In particular,  we have
$\chi (\omega _X )=1$.
So the claim will  follow from Proposition \ref{thetacohomology} if we show
that
$h^i(\omega_X\ot P)=0$ for all
$i>0$ and
$P\ne \OO _X$.

By the generic vanishing theorems (Theorem \ref{genvanish}),
the components of $V^i(\omega
_X)$ are torsion translates of subtori of $\Pic ^0(X)$ of codimension
at least $i$. Since $A$ is simple, it follows that
for all $i>0$ the sets  $V^i(\omega _X)$ consist of finitely
many torsion points
in $P\in \Pic ^0(X)$. Since $V^i(\omega_X)\subset V^1(\omega_X)$
for $i\ge 1$ by
Theorem
\ref{genvanish}, d), we have $V^i(\omega_X)=\{\OO_X\}$ for $i>0$ by Proposition
\ref{V1isolated}.\end{proof}
When $\Alb (X)$ is simple, it is also possible, using a result of Ein and
Lazarsfeld, to characterize abelian varieties in terms of $P_1(X)$.
\begin{lem}\label{V0}If $X$ is a smooth projective variety such that 
$P_1(X)>0$,
$q(X)=\dim (X)$
and $V^0(\omega _X)$ has dimension $0$, then
$X$ is birational to $\Alb(X)$.
\end{lem}
\begin{proof}
By \cite[ Prop. 2.2]{EL1}, $X \to \Alb (X)$ is surjective,
hence $X$ has maximal Albanese
dimension.

By an argument due to Ein and Lazarsfeld (cf. \cite[ Theorem 1.7]{CH1}),
$X$ is birational to $\Alb(X)$.
\end{proof}
\begin{cor}
Let $X$ be a nonsingular projective variety
of dimension $n$. If $X$ is of maximal Albanese dimension,
$\Alb (X)$ is simple, $P_1(X)=1$ and $q(X)=n$, then $X$ is
birational to $\Alb(X)$.
\end{cor}
\begin{proof} By assumption the Albanese map $\alb_X\colon X\to \Alb(X)$ is
generically finite and surjective.
By a result of Ueno (cf. \cite[3.4]{Mo}) one has $h^{i,0}(X)=
\binom {n}{i}$ for all $i\in [1,...,n]$. By Hodge symmetry and Serre
duality, one has that $h^{n-i}(\omega _X)=\binom {n}{i}$.
In particular $\chi (\omega _X)=0$. Since   $A$ is simple, $V^i(\ox)$  is 
$0-$dimensional for $i>0$.  Since
$\chi(\ox)=0$, then $V^0(\ox)$ is also $0-$dimensional 
and the claim follows from Lemma \ref{V0}.
\end{proof}
\section{The Albanese map}
Here we give a strengthening of a result
of Koll\`ar on the surjectivity of the Albanese map.
\begin{thm}\label{surjective}Let $X$ be a smooth projective variety.
If $P_2(X)=1$ or
$0<P_m(X)\leq 2m-3$ for some $m\geq 3$, then
$\alb _X\colon X\to \Alb (X)$ is surjective.
\end{thm}
\begin{proof}If $P_2(X)=1$ then this is \cite[Theorem 1]{CH3}.
If $P_3(X)=1$ or $0<P_m(X)\leq 2m-6$ for some $m\ge 4$, this is
\cite[Theorem 11.2]{Ko4}.
We will follow the proof of \cite[Theorem 11.2]{Ko4}.
Assume that $\alb _X\colon X\to \Alb (X)$ is not surjective.
By
\cite[10.9]{Ue}, up to replacing $X$ by an appropriate birational model,
there is a morphism $f\colon X
\to Z$ where
$Z$ is a smooth variety of general type of dimension $\ge 1$, such that the
Albanese map $a\colon Z\to S:=\Alb(Z)$ is birational onto its image.
We denote by $F_{X/Z}$
a general fiber of $f$.  In particular, notice that,
$\dim (S)\geq 2$.

\noindent{\bf Step 1.}\ \  {\em There exists an ample divisor  $L$ on
$Z$ such that, after replacing
$X$ by an appropriate birational model,
there exist an integer $r\gg 0$ and a divisor
$B\in |rm(K_X+(m-2)K_{X/Z})-f^*L)|$
such that $B$
has normal crossings support and
$$\lfloor\frac{ B|_{F_{X/Z}}}{rm}\rfloor \prec F$$
where $|mK_{F_{X/Z}}|=|H|+F$.}

We remark first of all that the linear system 
$|mK_{F_{X/Z}}|$ is nonempty, since
$P_m(X)>0$.  Let $L_0$ be an ample  $\QQ-$divisor on $Z$    such that
$K_Z-2L_0$ is big. By the proof of \cite[10.2]{Ko4}, after replacing $X$ by an
appropriate birational model, for $r$ sufficiently big and divisible  
there exists   a
divisor
$D\in |rm(m-1)K_{X/Z}+rmf^*L_0|$ such that the restriction 
of $D$ to the general fiber
$F_{X/Z}$ is equal to $\bar {H}+r(m-1)F$,
where  $|mK_{F_{X/Z}}|=|H|+F$ and $\bar {H}\in |r(m-1)H|$ is a
general smooth member. Let $G$ be a general member of
$|rm(K_Z-2L_0)|$ and set $L=rm L_0$. Let $B=f^*G+D\in |rm(K_X+(m-2)K_{X/Z}-f^*L)|$. After
replacing $X$ by an appropriate birational model, we may assume that $B$ has normal
crossings support. Then $B|_{F_{X/Z}}=
D|_{F_{X/Z}}=(\bar {H}+r(m-1)F).$
Therefore $$ \lfloor \frac {B|_{F_{X/Z}}}{rm}\rfloor=
\lfloor \frac {\bar{H}}{rm} + \frac {m-1}{m}F\rfloor\prec F.$$

\noindent {\bf Step 2.} {\em $\dim |(2K_X+(m-2)K_{X/Z})\ot  f^*P|\geq 1$
for all $P\in \Pic ^0(Z)$.}

Let $B$ be a divisor as in Step 1.
Define
$$M:=K_X+(m-2)K_{X/Z}-\lfloor \frac {B}{rm}\rfloor \equiv
\frac {f^*L}{mr}+\{ \frac {B}{rm}\}.$$
One has:
$$H^0(F_{X/Z},
\OO _{F_{X/Z}}(K_X+M))=
H^0(F_{X/Z},
\OO _{F_{X/Z}}(mK _{F_{X/Z}} -\lfloor \frac {B}{rm}\rfloor ))$$
\hskip2cm$$=H^0(F_{X/Z}, mK _{F_{X/Z}})>0.$$
In particular $f_*(K_X+M)\ne 0$, and hence, by Theorem \ref{kollar},   $h^0(f_
*(K_X+M)
\ot P)$ is a nonzero constant, independent of $P\in \Pic^0(Z)$.

We are going to show that this constant is $>1$. Indeed assume that it
is equal to $1$ and
consider
the Albanese map  $s\colon Z\to S:=\Alb(Z)$. By Lemma \ref{Ri} we have $h^i(S,
s_*f_*(K_X+M))=h^i(Z, f_*(K_X+M))$. If $i>0$,
then $h^i(Z, f_*(K_X+M))=0$ by Theorem
\ref{omegavan}, b),  and by
Proposition \ref{abelsupport} it follows that $s_*f_ *(K_X+M)$
is supported on an abelian subvariety of $S$.
However, by  Theorem \ref{omegavan}, a),
$f_ *(K_X+M)$
is a torsion free sheaf, and thus its support is $Z$.
Since $s$ is birational, the
support of
$s_*f_ *(K_X+M)$ is $s(Z)$, contradicting the fact  $Z$ is of general type.
This shows
that for every $P\in\Pic^0(Z)$ we have
$$h^0((2K_X+(m-2)K_{X/Z})\ot P)\geq  h^0(f_ *(K_X+M)\ot P)\geq 2.$$

\noindent {\bf Step 3.}

 We have  (cf. \cite[11.2]{Ko4}) that  $P_1(Z)\geq \dim (Z)+1$,
$P_r(Z)\geq 2r-1$ for all $r\geq 2$, and if $\dim (Z)\geq 2$
then $P_r(Z)\geq 2r$.
 It follows that by
Lemma \ref{Pm},
$h^0(rK_Z+P)\geq 2r-1$ for all $ r \geq 2$ and $P\in \Pic ^0(Z)$, and if
$\dim (Z)\geq 2$
then $h^0(rK_Z+P)\geq 2r$.

We now apply Lemma \ref{claimA} with $T:=f^*\Pic^0(Z)$,
$L= 2K_X+(m-2)K_{X/Z}$, $M=(m-2)f^*K_Z$, so that $L\ot M=mK_X$.
Recalling that $\dim\Pic^0(Z)\ge 2$, for $m\ge 4$ we get $\dim|mK_X|\ge
2(m-2)-2+1+2=2m-3$, the required contradiction.
For $m=3$ and $\dim (Z)=1$, we have $h^0(K_Z-P)\geq 1$ for every
$P\in \Pic^0(Z)$
and Lemma \ref{claimA} gives again $\dim |3K_X|\ge 3=2m-3$.
Finally, for  $m=3$ and  $\dim (Z)\geq 2$, we have $h^0(K_Z)\geq 3$ and the
contradiction follows by considering   the morphism of linear series
$$|2K_X+K_{X/Z}|\times f^*|K_Z|\to |3K_X|.$$
\end{proof}
It is natural to ask whether these bounds are optimal. For example,
does there exists an $X$ with $P_3(X)=4$, such that $\alb _X$ is not
surjective?

\section{Varieties with $P_3(X)=2$ and $q(X)=\dim (X)$}\label{P3=2}
In this section we give an explicit birational description of the varieties  
$X$ with
$P_3(X)=2$ and
$q(X)=\dim (X)$. A key step  of the proof is the use of 
Theorem \ref{surjective}.  The
precise statement of the result  is as follows:
\begin{thm} Let $X$ be a smooth projective variety.
\noindent Then
$P_3(X)=2$ and $q(X)=\dim (X)$ iff:
\begin{itemize}
\item[a)] there is a surjective map $\phi\colon \Alb(X)\to E$, 
where $E$ is a curve of
genus
$1$;

\item[b)]$\alb _X\colon X\to\Alb (X)$
is birational to a smooth double cover of $\Alb (X)$ defined by
$Spec (\OO _{\Alb (X)}\oplus P\otimes\phi^*L\inv)$, 
where $L$ is a line bundle of $E$ of
degree $1$ and
$P\in
\Pic ^0(X)\setminus \phi^*\Pic^0(E)$ is
$2$-torsion. The branch locus of the double cover is the union of 
two distinct fibers of
$\phi$.
\end{itemize}
\end{thm}\begin{proof}
If $X$ is the double cover described in the statement, 
then  the  standard formulas for
double covers give $P_1(X)=1$ and, for $m\ge 2$,  
$P_m(X)=m$ if $m$ is even and $P_m(X)=m-1$
if
$m$ is odd.

Assume now that $P_3(X)=2$ and $q(X)=\dim (X)$.
By Theorem \ref{surjective},
$X\to
\Alb (X)$ is surjective and thus $X$ has maximal Albanese dimension.
Let $f\colon X\to Y$ denote the Iitaka fibration of $X$. We have a commutative
diagram:
$$
\CD
X &@>\alb _X>>&\Alb (X) \\
@V{f}VV &  &  @VV{f_*}V \\
Y & @>\alb _Y>>&\Alb (Y).\\
\endCD
$$
Moreover, by Proposition \ref{albanese}, $K:= ker f_*$ is connected
and
there exists an abelian variety $P$ isogenous to $K$ and
birational to $F:=F_{X/Y}$.
Let $$G:=ker \left( \Pic ^0(X)\to \Pic ^0(F)\right).$$
Then $G$ is the union of finitely many translates of $\Pic ^0(Y)$
corresponding to the finite group $$G/\Pic ^0(Y)\cong ker \left(
\Pic ^0(K)\to \Pic ^0(F)\right).$$
For any line bundle $L$ on $X$, define (cf. \S 2)
$$V^i(L):=\{
P\in \Pic ^0(X) | h^i(L\ot P)\ne 0\} .$$
Since $(\omega _X^m\ot P)|F=\omega^m _F\ot P$, it follows that
$V^0(\omega _X^m)\subseteq G$ for every $m\ge 0$. By \cite{CH2} Lemma 2.2,
one sees
that for all $P\in G\setminus\Pic ^0(Y)$,
the dimension of $V^0(\omega _X)\cap (P+\Pic ^0(Y))$ is $\ge 1$.
Moreover, by Corollary \ref{Vzero}
the only possible $0$-dimensional component of $V^0(\omega _X)$
is the origin.
The proof is divided into several steps.

\noindent{\bf Step 1.} {\em $V^0(\omega _X ^2)= G$.}

Since $h^0(\omega_X)>0$,  we have $V^0(\omega_X)\subset V^0(\omega _X ^2)$.
By Proposition  \ref{Pm}, we have $h^0(2K_X+ P)= P_2(X)>0$ for
every
$P\in
\Pic^0(Y)$, namely $V^0(\omega_X^2)\cap \Pic^0(Y)=\Pic^0(Y)$. Consider an
irreducible component
$S=Q+\Pic^0(Y)$ of
$G$,  where $Q\in \Pic^0(X)\setminus \Pic^0(Y)$ is torsion. 
Then by \cite[Lemma 2.2]{CH2}
$V^0(\omega _X)\cap S$ has positive dimension, in particular it is nonempty.
Thus $V^0(\omega _X^2)\cap S$ is also nonempty and   Proposition
\ref{Pm} implies again $V^0(\omega _X^2)\cap S=S$.
\medskip

\noindent{\bf Step 2.} {\em If $T$ is an irreducible component of $V^0(\omega
_X)$, then $\dim (T)\leq 1$.}

Since $T\subset G$,
we have also $-T\subset G$, and hence by Step 1, 
$h^0(2K_X-Q)>0$ for all $Q\in T$.
Applying Lemma \ref{claimA} with $L=K_X$ and $M=2K_X$ we get $\dim
(T)\le \dim|3K_X|=1$.
\medskip

\noindent{\bf Step 3.} {\em For every $P\in V^0(\ox)$ one has $h^0(K_X+P)=1$.}

 Assume that there is $P\in V^0(\ox)$ such that  $h^0(K_X+ P)\geq 2$.
Then we have
$h^0(2K_X+2P)\geq 3$ and thus,
by Proposition \ref{Pm}, we have  $h^0(2K_X+2P+R) \geq 3$ for all 
$R\in \Pic ^0(Y)$. Since
$P\in G$,  then $-2P\in G$. Let
$P_0$ be a point  of $V^0(\omega _X)\cap -2P +\Pic ^0(Y)$ 
(such a point exists by \cite[Lemma
2.2]{CH2}). Then from
$h^0(2K_X-P_0) \geq 3$ and $h^0(K_X+P_0)>0$  it follows immediately
$h^0(3K_X)\ge 3$, against the assumptions.
\medskip

\noindent{\bf Step 4.} {\em For any $P\in G\setminus\Pic ^0(Y)$, 
$h^0(2K_X+P)\leq 1$.}

Let $P\in G\setminus\Pic^0(Y)$ be such that
 $h^0(2K_X+P)\ge 2$. By Proposition \ref{Pm}, one has $h^0(2K_X+R+P) \geq 2$ 
for all
$R\in \Pic ^0(Y)$.
Since $P\in G\setminus\Pic^0(Y)$ also $-P \in G\setminus\Pic^0(Y)$. 
Let $T$ be an irreducible
component of $V^0(\omega _X)\cap -P +\Pic ^0(Y)$. Then by Step 1 and 
Proposition \ref{Pm}, one has
$h^0(2K_X-Q) \geq 2$ for all $Q\in T$.
Since $T$ has positive dimension by \cite[Lemma 2.2]{CH2}, 
Lemma \ref{claimA} gives
$\dim|3K_X|\ge 2$, a contradiction.
\medskip

\noindent{\bf Step 5.} {\em $V^0(\omega _X)\cap \Pic ^0(Y)=\{ \OO _X\}$.}

Recall that by Corollary \ref{Vzero} $\{\OO_X\}$ is the only possible
$0-$dimensional component of $V^0(\ox)$.
Let
$T$ be a positive dimensional component of
$V^0(\omega _X)\cap \Pic ^0(Y)$. Then by Theorem \ref{genvanish}, b),  
$T=Q+T_0$, where $Q\in
\Pic^{\tau}(Y)$ and
$T_0$ is an abelian subvariety of $\Pic^0(Y)$. For every $P\in T_0$, we have
$h^0(K_X+Q+P)>0$ and $h^0(K_X+Q-P)>0$. Thus   Lemma
\ref{claimA} gives $h^0(2K_X+2Q)\ge 2$. By Proposition \ref{Pm} we have 
$h^0(2K_X+P)\ge
2$ for every $P\in \Pic^0(Y)$. In particular, we have $h^0(2K_X-P)\ge 2$ 
and $h^0(K_X+P)\ge 1$
for all $P\in T$.  Applying Lemma
\ref{claimA} again gives
$\dim|3K_X|\ge 2$, a contradiction.
\medskip

\noindent{\bf Step 6.} {\em $Y\to\Alb (Y)$ is birational.}

Since $X$ is of Albanese general type, the same is true for
$Y$ by Proposition \ref{albanese}. 
So if $\kappa(Y)=0$ then the claim is true by
\cite{Ka}. Assume
$\kappa(Y)>0$.    By the commutativity of the diagram 
at the beginning of the proof,  the map  $Y\to
\Alb(Y)$ is surjective,  since
$X\to \Alb(X)$ is surjective. By a result
of Ein and Lazarsfeld (see also
\cite[ Theorem 2.3]{CH2}) there is an irreducible component
$T\subset V^0(Y,\omega _Y)$ of positive dimension. 
The divisor $K_{X/Y}$ is effective since
$X$ has maximal Albanese dimension, 
hence $\dim \left( V^0(X,\omega _X)\cap \Pic
^0(Y)
\right)
\geq 1$, contradicting  Step 5.
\medskip

Denote by  $\pi \colon X\to \Alb (Y)$ the composition of the Albanese map $X\to
\Alb(X)$ and of $f_*\colon \Alb(X)\to\Alb(Y)$.

\noindent{\bf Step 7.} {\em For any $P\in G\setminus\Pic ^0(Y)$, there exists
a principal polarization $N$ on $\Alb (Y)$ and an effective divisor $D$
on $X$ such that
$$ |2K_X +P|=\pi ^*|N|+D.$$}

 See Step 3 in the proof of \cite[Theorem 2.4]{CH1}.
\medskip

Let $T\subset \Pic ^0(X)$
be any $1-$dimensional irreducible component of $V^0(\omega _X)$. 
We recall that by
Theorem
\ref{genvanish}, b),
$T$ is a translate of an abelian subvariety
$\bar{T}$.  Let $\Alb (X)\to E:=
\Pic^0(\bar{T})$ be the dual map of abelian varieties
and $\pi_E \colon X\to E$  the induced morphism.
By Lemma \ref{lemmaD}, there exists a divisor 
$D\prec R:=Ram (\alb _X)=K_X$, vertical
with respect to $\pi_E$, such that for all $P\in T$, 
$F:=R-D$ is a fixed divisor
of each of the linear series $|K_X+P|$.
In other words, we have for all $P\in T$
$$|K_X+P|=F+|V_P|.$$
The divisors in $|V_P|$ are vertical
with respect to $\pi_E$. To see this, 
notice that $(K_X-F)|_{F_{X/E}}=D|_{F_{X/E}}=\OO _{
F_{X/E}}$ so $(K_X-F+P)|_{F_{X/E}}=(V_P)|_{F_{X/E}}=\OO _{
F_{X/E}}$.

Since $\bar{T}\subset \Pic^0(Y)\subset \Pic^0(X)$,  the map $\pi_E$
factors through the  map $\Alb(Y)\to E$, which has connected fibers since it is
dual to an inclusion of tori. The map
$X\to \Alb(Y)$ has connected fibers by Step 6, hence
$\pi_E$ has connected fibers.
\medskip

\noindent{\bf Step 8.} {\em  For general $P\in T$, there exists
a line bundle  $L_P$  of degree $1$ on
$E$ such that $\pi_E^*L_P\prec K_X+P$.}

The divisors $V_P$ are nonempty and vary with $P\in T$. Since $V_P$ 
is vertical, for $P\in
T$ general it contains a smooth fiber of $\pi(e)$, where $e\in E$. 
So we may set
$L_P:=\OO_{E}(e)$.
\medskip

\noindent{\bf Step 9.} {\em $\kappa (X)=1$.}

Let $T$ be a component of $V^0(\ox)$ of positive dimension. By Step 2,
$T$ has dimension 1. As before, we let $E=\Pic^0(\bar{T})$ 
and we denote by $\pi_E\colon
X\to E$ the corresponding map.  Let
$P\in T$ be a general point. By Step 8, there is a line  bundle $L_P$  of
degree 1 on
$E$ such that $\pi_E^*{L_P}\prec K_X+P$, and by Step 7 there is a
 principal polarization $N$  on
$\Alb(Y)$ such that $\pi^*N\prec 2K_X-P$. We denote by $L$ the pull-back of
${L_P}$ to
$\Alb(Y)$. There is an inclusion $\pi^*|L\ot N|\to |3K_X|$, 
and thus the dimension of $|L\ot
N|$ is $\le 1$. On the other hand, for every $Q\in \bar{T}=\Pic^0(E)\subset
\Pic^0(Y)$ we have
$h^0(\Alb(Y),L+ Q)=1$ and $h^0(\Alb(Y),N-Q)=1$. So Lemma
\ref{claimA} gives $\dim|L\ot
N|\ge  1$. It follows that $\dim|L\ot N|=1$, and the proof of 
Lemma \ref{claimA} shows that
for  every divisor D in $|L\ot N|$ there is $Q\in \Pic^0(E)$ 
such that $D$ can be written as
$D'_Q+D''_Q$, where
$D'_Q\in |L+Q|$ and $D''_Q\in|N-Q|$. Let $|M|$ be the moving part 
of $|L\ot N|$. The line
bundle $M$ is positive semidefinite, so that  there exist a 
quotient $\bar{A}$ of
$\Alb(Y)$ and an ample line bundle $\bar{M}$ on $\bar{A}$ such that $M$ 
is the pull-back of
$\bar{M}$ to
$\Alb(Y)$ and $h^0(M)=h^0(\bar{M})$ (cf. \cite{LB}, \S 3.3). 
Notice that since the continuous
system
$\{D'_Q\}$ is base point free, every divisor of 
$|M|$ can be written as $D'_Q+R_Q$ for
suitable
$Q\in\Pic^0(E)$. Moreover, $R_Q>0$, since 
$0=\dim|L_Q|<1=\dim|M|$. So the general divisor
of $|M|$ is reducible, and the same is true 
for the general divisor of $|\bar{M}|$ on
$\bar{A}$. Since $\bar{M}$ is ample and has no fixed part, 
by Theorem 4.5 of \cite{LB} this
can only happen if $\bar{A}$ has dimension 1.  
Since $L\prec M$, it follows that
$\bar{A}=E$, namely the moving part of $|L\ot N|$, 
and thus also of $|3K_X|$, is a
pull-back from $E$. This condition determines 
the map $\Alb(Y)\to E$ uniquely and, taking
duals, it determines uniquely $\bar{T}\subset\Pic^0(Y)$. 
So we have shown that all the
positive dimensional components of $V^0(\ox)$ are 
translate of the same abelian subvariety
$\bar{T}\subset \Pic^0(Y)$ of dimension 1. Since
by \cite[Theorem 2.3]{CH2},
$\Pic ^0(Y)$ is generated by the sum of these translates, 
it follows that $\Alb(Y)$ has
dimension 1, namely
$\kappa(X)=1$.
\medskip

\noindent{\bf Step 10.} {\em $V^0(\ox)=G\setminus \Pic^0(Y)\cup\{\OO_X\}$.}

By Step 6 and Step 9, $Y$ is a smooth curve of genus 1. 
As we have remarked at the beginning
of the proof, $V^0(\ox)$ intersects each component of 
$G\setminus\Pic^0(Y)$ in a set of
positive dimension. This shows that $V^0(\ox)$ contains 
$G\setminus\Pic^0(Y)$. The statement
now follows from Step 5.
\medskip

\noindent{\bf Step 11.} {\em  $\alb _X$ is of degree 2.}

Since $X\to\Alb(Y)$ and $\Alb (X)\to\Alb(Y)$
have connected fibers, then
$$d=\deg \left( X\to \Alb(X) \right)=  \deg \left( F_{X/Y}\to K \right)$$
is just the cardinality of $G/\Pic ^0(Y)$.
If $d\geq 3$, then there exist elements 
$P_1,P_2,P_3 \in V^0(\omega _X)\setminus
\Pic ^0(Y)$ such that $P_1+P_2+P_3=\OO _X$.
Since
$h^0(K_X + P_i+ Q)=1$ for all $Q\in \Pic ^0(Y)$ by Step 10,
it follows applying Lemma \ref{claimA}
that $h^0(2K_X+ P_1+P_2 +Q)\geq 2$ for all $Q\in \Pic ^0(Y)$
and similarly that
$h^0(3K_X)=h^0(3K_X+P_1+P_2+P_3)\geq 3$, a contradiction.
\medskip

\noindent{\bf Step 12.} 
{\em ${\alb _X}_*\omega _X=\OO _A \oplus (f_*)^* L\ot P$
where $L$ is a principal polarization on $\Alb(Y)$ and
$P^{\ot 2}=\OO _X$ but $P\notin \Pic^0(Y)$.}

By Step 10 and Step 11, we have
$$V^0(\omega _X)=\{ \OO _X\} \cup( P+\Pic ^0(Y))$$
for an appropriate 2-torsion element $P\in \Pic ^0(X)\setminus \Pic^0(Y)$.
The sheaf $L:=f _* (\omega _X \ot P)$ is torsion 
free by Theorem \ref{omegavan}, a) and
it has rank 1 since the general fiber $F_{X/Y}$  
of $f$ is birational to an abelian
variety and the restriction of  $P$ to $F_{X/Y}$ 
is trivial. Since $Y$ is a curve, $L$
is actually a line bundle. In addition, we have
$h^0(L\ot Q)=1$ for all $Q\in \Pic ^0(Y)$. Therefore
$L$ is a principal polarization.
Similarly the sheaf $M:=f _* (\omega _X )$ is
a line bundle
such that $h^0(M)=1$ and $h^0(M\ot Q)=0$ for all $Q\in \Pic ^0(Y)\setminus\{\OO
_Y\}$, therefore $M=\OO _Y$.
Let $f_*\colon \Alb (X)\to \Alb (Y)=Y$ and let $\tilde{L}:=(f_*)^*L$. There are
inclusions of sheaves
$$\OO_{\Alb(X)}=(f_*) ^*(f_*) _* ({\alb _X}_* \omega _X)
\to {\alb _X}_* \omega _X,$$
$$\tilde{L}=(f_*) ^*(f_*) _* 
({\alb _X}_* \omega _X\otimes P)\to {\alb _X}_* (\omega _X
\otimes P).$$
There is a corresponding map
$$\psi \colon  \OO _{\Alb(X)}\oplus \tilde{L} \ot P \to {\alb _X}_*\omega _X.$$
Restricting to a generic fiber
$K$ of $\Alb(X)\to \Alb (Y)=Y$, one has that
the above map is given by
$$\OO _K\oplus P|K \to ({\alb _X}_*\omega _X)|K.$$
Let $F$ be a general fiber of $X\to Y$,
$a:={\alb _X}|_F\colon F\to K$. The above map of sheaves
is equivalent to the isomorphism
$$\OO _K\oplus P|K\cong a_*\omega _F .$$
It follows that $\psi $ is an inclusion of sheaves.

We wish to show that $\psi$ is an isomorphism.
To this end, by Corollary 2.2,
it suffices to show that $\psi$ induces isomorphisms
of cohomology groups
$$
H^i\left( (\OO _{\Alb(X)} \oplus \tilde{L}\otimes P)\otimes Q\right) \to
H^i\left({\alb _X}_* (\omega _X) \ot Q\right)$$
for all $i\ge 0$ and $Q\in \Pic ^0(X)$.

For $i=0$ the map above is injective and it is enough to check that the 
two vector spaces
have the same dimension. By Step 3 and by the
description of $V^0(\ox)$ that we have given, it follows that
for $i=0$ and any $Q\in \Pic ^0(X)$,
we have isomorphisms in cohomology.  If $i\geq 1$ and $Q$ is not
in $\{\OO _X\}\cup (\Pic^0(Y)+P)$,
then both vector spaces are $0$ (cf. Theorem
\ref{genvanish}, d)).  We will prove that the above isomorphism
holds for all
$ Q\in P+\Pic ^0(E)$. (The proof for $Q=\OO _X$ proceeds analogously
but is easier).
For such a choice of $Q$ one has:
$$H^i\left(\tilde{L}\ot P \otimes Q\right) =
H^i\left( (\OO _A \oplus  \tilde{L}\ot P)\otimes Q\right).$$
We observe that $\chi(\ox)=0$ by Theorem \ref{genvanish}, d), since $X$ 
has maximal Albanese
dimension and we have seen that
$V^0(\ox)$ is a proper subset of
$\Pic^0(X)$.
Let $W\subset H^1(\OO_X)$ be a subspace complementary to the tangent space 
to $\Pic^0(Y)$.
The assumptions of  Proposition  \ref{gencoho} are satisfied for all $Q\in
P+\Pic^0(E)$, hence  for all
$0\leq j\leq
\dim (X)$ and
$Q\in P+\Pic ^0(E)$ there are isomorphisms
$$H^j(\omega _X\ot Q)\cong H^0(\omega _X\ot Q)\ot
\wedge ^{j} W.$$ induced by cup product.
The required  isomorphism  is given by the following commutative diagram
 $$
\CD
H^0(\tilde{L} \ot P\ot  Q) \ot \wedge ^{i}  {W}
& @>{\cong}>> &
H^0({\alb _X}_* \omega _X \ot Q)\ot \wedge ^{i}{W} \\
@V{\cong}VV & & @VV{\cong}V & \\
H^i(\tilde{L}\ot P\ot  Q)& @>>> &
H^i(\alb _* \omega _X \ot Q) \\
\endCD
$$
\medskip

\noindent{\bf Step 13.} {\em  Conclusion of the proof.}

If $X\overset{\epsi}{\to} Z\overset{g}{\to}\Alb(X)$ is the Stein 
factorization of  $\alb_X$,
then by  Proposition \ref{OXfree}, $g$ is flat and 
$g_*\OO_Z=\OO_{\Alb(X)}\oplus
\tilde{L}\inv
\ot P$. The branch locus of $g$ is reduced, since $Z$ is normal, 
and it is linearly
equivalent to $2\tilde{L}$, thus it consists of two fibers of $f_*$.
\end{proof}

\bigskip
\bigskip

\begin{minipage}{13cm}
\parbox[t]{5.5cm}{Christopher D. Hacon\\
Department of Mathematics\\
Sproul Hall 2208\\
University of California\\
Riverside, CA 92521-013 USA\\
hacon@math.ucr.edu}
\hfill
\parbox[t]{5.5cm}{Rita Pardini\\
Dipartimento di Matematica\\
Universit\`a di Pisa\\
Via Buonarroti, 2\\
56127 Pisa, ITALY\\
pardini@dm.unipi.it}
\end{minipage}

\end{document}